\documentclass[12pt]{article}
\usepackage[top=1in, bottom=1in, left=1in, right=1in, a4paper]{geometry}

\usepackage[english]{babel}

\usepackage{amssymb,amsmath,amsthm}
\usepackage{graphicx,cite,color}
\usepackage{longtable,pdflscape,booktabs,caption,multicol}
\usepackage[colorlinks=true,citecolor=black,linkcolor=black,urlcolor=blue]{hyperref}
\usepackage{enumitem}
\usepackage{mathtools, bm}

\theoremstyle{plain}
\newtheorem{theorem}{Theorem}
\newtheorem{lemma}[theorem]{Lemma}

\newtheorem{proposition}[theorem]{Proposition}

\newtheorem{problem}[theorem]{Problem}

\theoremstyle{definition}
\newtheorem{definition}[theorem]{Definition}

\theoremstyle{remark}
\newtheorem*{remark}{Remark}

\usepackage{float, tikz, subfig}
\usetikzlibrary{automata}
\usetikzlibrary{decorations.pathmorphing}

\usepackage{listings}
\definecolor{mauve}{rgb}{0.58,0,0.82}
\definecolor{dkgreen}{rgb}{0,0.6,0}
\lstdefinestyle{pitonche} {
    language = Python,
    basicstyle = footnotesizettfamily,
    showspaces = false,
    showstringspaces = false,
    breakautoindent = true,
    flexiblecolumns = true,
    keepspaces = true,
    stepnumber = 1,
    xleftmargin = 0pt
}
\usepackage{authblk}

\title{Characterization of transmission irregular\\ starlike and double starlike trees}

\author[1,2]{Ivan Damnjanović}
\affil[1]{Faculty of Electronic Engineering, University of Niš, Niš, Serbia}
\affil[2]{Diffine LLC, San Diego, California, USA}
\affil[ ]{{\tt ivan.damnjanovic@elfak.ni.ac.rs}}

\date{}

\begin{document}

\maketitle

\begin{abstract}
The transmission of a vertex in a connected graph is the sum of its distances to all the other vertices. A graph is transmission irregular (TI) when all of its vertices have mutually distinct transmissions. In an earlier paper, Al-Yakoob and Stevanović [Appl.\ Math.\ Comput.\ {\bf 380} (2020), 125257] gave the full characterization of TI starlike trees with three branches. Here, we improve these results by using a different approach to provide the complete characterization of all TI starlike trees. Moreover, we find the precise conditions under which a double starlike tree is TI. Finally, we implement the aforementioned conditions in order to find several infinite families of TI starlike trees and TI double starlike trees.
\end{abstract}


\section{Introduction}

Let $G$ be a simple connected graph with the vertex set $V(G)$. The \emph{degree} $d_G(u)$ of the vertex $u \in V(G)$ is the total number of vertices in $V(G)$ adjacent to $u$. The \emph{distance} $d_G(u, v)$ between the vertices $u, v \in V(G)$ is the length of the shortest walk between $u$ and $v$ in $G$. The \emph{transmission}
\[
    \mathrm{Tr}_G(u) = \sum_{v \in V(G)} d_G(u, v)
\]
of the vertex $u \in V(G)$ is the sum of its distances to all the other vertices of $G$. For the sake of brevity, we will drop the subscripts from $d_G$ and $\mathrm{Tr}_G$ whenever the graph is clear from the context. A graph is \emph{transmission irregular}, or TI for short, if all of its vertices have distinct transmissions.

Alizadeh and Klavžar \cite{AlKl2018} have shown that almost no graph is transmission irregular. This is a direct consequence of the fact that almost all graphs have diameter two \cite{BlHa1979}, bearing in mind that for any such graph $G$ we have
\[
    \mathrm{Tr}_G(u) = 2\left( |V(G)| - 1 \right) - d_G(u) \qquad \left( u \in V(G) \right) ,
\]
while it is guaranteed that at least two vertices have the same degree. Moreover, Alizadeh and Klavžar \cite{AlKl2018} have also investigated the TI property of starlike trees with three branches.

\begin{definition}
    A tree $T$ is \emph{starlike} if it contains a single vertex $u$ of degree greater than two. If $d_T(u) = k$ and the pendent paths starting at $u$ have lengths $A_1, A_2, \ldots, A_k \in \mathbb{N}$, then we will denote such a starlike tree by $\mathrm{S}(A_1, A_2, \ldots, A_k)$.
\end{definition}

The following theorem provides the full characterization of starlike trees bearing the form $\mathrm{S}(1, A_2, A_3)$.

\begin{theorem}[\hspace{1sp}{\cite[Theorem 2.2]{AlKl2018}}]\label{AlKl_th}
    If $1 < A_2 < A_3$, then $\mathrm{S}(1, A_2, A_3)$ is TI if and only if $A_3 = A_2 + 1$ and $A_2 \notin \left\{ \frac{k^2-1}{2}, \frac{k^2-2}{2} \right\}$ for some $k \ge 3$.
\end{theorem}

Theorem \ref{AlKl_th} was further improved by Al-Yakoob and Stevanović \cite[Theorem 2]{AlYaSt2020}, who disclosed the full characterization of all starlike trees $\mathrm{S}(A_1, A_2, A_3)$ with three branches. Besides inspecting the TI property of given classes of graphs, many researchers have found interest in searching for various families of TI graphs (see, for example, \cite{Dobrynin2019a, Dobrynin2019b, Dobrynin2019c, AlYaSt2020, AlYaSt2022a, AlYaSt2022b, BeDo2021, BeDo2022, DoSha2020, XuKl2021, XuTiKl2023}). Xu et al.\ \cite{XuTiKl2023} have recently used an approach based on double starlike trees to investigate the TI property of chemical trees.

\begin{definition}
    A graph is \emph{chemical} if its maximum vertex degree is not greater than four.
\end{definition}
\begin{definition}
    A tree $T$ is \emph{double starlike} if it contains exactly two vertices $u$ and $v$ whose degree is greater than two. If $d_T(u) = k + 1, \, d_T(v) = m + 1, \, d_T(u, v) = C$ and the pendent paths starting at $u$ have lengths $A_1, A_2, \ldots, A_k$, while the pendent paths starting at $v$ have lengths $B_1, B_2, \ldots, B_m$, we will signify such a double starlike tree as $\mathrm{DS}(C; A_1, A_2, \ldots, A_k; B_1, B_2, \ldots, B_m)$.
\end{definition}

Without loss of generality, we will always assume that $\sum_{i=1}^k A_i \ge \sum_{i=1}^m B_i$ for each double starlike tree. The next problem was subsequently posed by Xu et al.

\begin{problem}[\hspace{1sp}{\cite[Problem 5.1]{XuTiKl2023}}]\label{XuTiKl_prob}
    Does there exist a TI chemical graph of every even order?
\end{problem}

Motivated by the aforementioned results, we begin by providing the complete characterization of all TI graphs among starlike and double starlike trees, thereby extending the theorem previously derived by Al-Yakoob and Stevanović \cite[Theorem 2]{AlYaSt2020}. Before stating our main results, let $g \colon \mathbb{R} \times (0, +\infty) \to \overline{\mathbb{R}}$ be the auxiliary function defined via
\[
    g(C_1, C_2) = \begin{cases}
        \frac{1}{2} \left( C_2 + \sqrt{C_2^2 - 4 C_1} \right), & C_2^2 - 4 C_1 \ge 0,\\
        -\infty, & C_2^2 - 4C_1 < 0 .
    \end{cases}
\]

The complete characterization of TI starlike and double starlike trees is disclosed within the next two theorems.

\begin{theorem}\label{main_th_1}
    For the starlike tree $\mathrm{S}(A_1, A_2, \ldots, A_k)$, let
    \[
        \alpha_{ij} = (n - 1 - A_i - A_j)(A_j - A_i) \qquad (1 \le i, j \le k).
    \]
    The given tree is transmission irregular if and only if the next conditions hold:
    \begin{enumerate}[label=\textbf{(\roman*)}]
        \item the values $A_1, A_2, \ldots, A_k$ are mutually distinct;
        \item $A_i < \frac{n}{2}$ for each $1 \le i \le k$;
        \item for all the $1 \le i < j \le k$, the number $\alpha_{ij}$ does not contain a positive divisor $p \in \mathbb{N}, \, p > \sqrt{|\alpha_{ij}|}$ such that the following conditions are satisfied:
        \begin{enumerate}[label=\textbf{(\alph*)}]
            \item $p$ and $\frac{\alpha_{ij}}{p}$ are of the same parity if and only if $2 \nmid n$;
            \item $g(\alpha_{ij}, n + 1 - 2A_i) \le p \le g(\alpha_{ij}, n-1)$;
            \item $g(\alpha_{ji}, n + 1 - 2A_j) \le p \le g(\alpha_{ji}, n-1)$.
        \end{enumerate}
    \end{enumerate}
\end{theorem}

\begin{theorem}\label{main_th_2}
    For the double starlike tree $\mathrm{DS}(C; A_1, A_2, \ldots, A_k; B_1, B_2, \ldots, B_m)$, let $A_* = \sum_{i=1}^k A_i$ and $B_* = \sum_{i=1}^m B_i$. Furthermore, let
    \begin{alignat*}{2}
        \alpha_{ij} &= (n - 1 - A_i - A_j)(A_j - A_i) \qquad & (1 \le i, j \le k) ,\\
        \beta_{ij} &= (n - 1 - B_i - B_j)(B_j - B_i) \qquad & (1 \le i, j \le m) ,\\
        \gamma_{i} &= (n - 1 - A_i - A_*)(A_* - A_i) \qquad & (1 \le i \le k) ,\\
        \delta_{ij} &= (n - 1 - A_i - B_j)(B_j - A_i) + C(A_* - B_*) \qquad & (1 \le i \le k, 1 \le j \le m) .
    \end{alignat*}
    The given tree is transmission irregular if and only if the next conditions hold:
    \begin{enumerate}[label=\textbf{(\roman*)}]
        \item the values $A_1, A_2, \ldots, A_k$ are mutually distinct and the values $B_1, B_2, \ldots, B_m$ are also mutually distinct;
        \item $A_i < \frac{n}{2}$ for each $1 \le i \le k$, $B_i < \frac{n}{2}$ for each $1 \le i \le m$, as well as $1 + A_* > \frac{n}{2}$;
        \item for all the $1 \le i < j \le k$, the number $\alpha_{ij}$ does not contain a positive divisor $p \in \mathbb{N}, \, p > \sqrt{|\alpha_{ij}|}$ such that the following conditions are satisfied:
        \begin{enumerate}[label=\textbf{(\alph*)}]
            \item $p$ and $\frac{\alpha_{ij}}{p}$ are of the same parity if and only if $2 \nmid n$;
            \item $g(\alpha_{ij}, n + 1 - 2A_i) \le p \le g(\alpha_{ij}, n-1)$;
            \item $g(\alpha_{ji}, n + 1 - 2A_j) \le p \le g(\alpha_{ji}, n-1)$;
        \end{enumerate}
        \item for all the $1 \le i < j \le m$, the number $\beta_{ij}$ does not contain a positive divisor $p \in \mathbb{N}, \, p > \sqrt{|\beta_{ij}|}$ such that the following conditions are satisfied:
        \begin{enumerate}[label=\textbf{(\alph*)}]
            \item $p$ and $\frac{\beta_{ij}}{p}$ are of the same parity if and only if $2 \nmid n$;
            \item $g(\beta_{ij}, n + 1 - 2B_i) \le p \le g(\beta_{ij}, n-1)$;
            \item $g(\beta_{ji}, n + 1 - 2B_j) \le p \le g(\beta_{ji}, n-1)$;
        \end{enumerate}
        \item for all the $1 \le i \le k$, the number $\gamma_{i}$ does not contain a positive divisor $p \in \mathbb{N}, \, p > \sqrt{|\gamma_{i}|}$ such that the following conditions are satisfied:
        \begin{enumerate}[label=\textbf{(\alph*)}]
            \item $p$ and $\frac{\gamma_{i}}{p}$ are of the same parity if and only if $2 \nmid n$;
            \item $g(\gamma_{i}, n + 1 - 2A_i) \le p \le g(\gamma_{i}, n-1)$;
            \item $g(-\gamma_{i}, 2A_* - n + 3) \le p \le g(-\gamma_{i}, 2A_* + 2C - n + 1)$;
        \end{enumerate}
        \item for all the $1 \le i \le k, 1 \le j \le m$, the number $\delta_{ij}$ does not contain a positive divisor $p \in \mathbb{N}, \, p > \sqrt{|\delta_{ij}|}$ such that the following conditions are satisfied:
        \begin{enumerate}[label=\textbf{(\alph*)}]
            \item $p$ and $\frac{\delta_{ij}}{p}$ are of the same parity if and only if $2 \nmid n$;
            \item $g(\delta_{ij}, n + 1 - 2A_i) \le p \le g(\delta_{ij}, n-1)$;
            \item $g(-\delta_{ij}, n + 1 - 2B_j) \le p \le g(-\delta_{ij}, n-1)$.
        \end{enumerate}
    \end{enumerate}
\end{theorem}

We further apply Theorems \ref{main_th_1} and \ref{main_th_2} in order to yield various families of TI trees providing a partial answer to Problem \ref{XuTiKl_prob}. For convenience, let an \emph{X-tree} $\mathrm{X}(A_1, A_2, A_3, A_4)$ be a starlike tree with four branches and let an \emph{H-tree} $\mathrm{H}(C; A_1, A_2; B_1, B_2)$ be a double starlike tree with two vertices of degree three. We disclose the following four infinite families of X-trees composed solely of TI chemical trees of even order.

\begin{theorem}\label{x_tree_th}
    For any $t \in \mathbb{N}_0$, each of the four X-trees
    \begin{enumerate}[label=\textbf{(\roman*)}]
        \item $\mathrm{X}(7 + 12t, 6 + 12t, 3 + 12t, 1)$,
        \item $\mathrm{X}(9 + 12t, 8 + 12t, 5 + 12t, 3)$,
        \item $\mathrm{X}(10 + 12t, 9 + 12t, 6 + 12t, 2)$,
        \item $\mathrm{X}(11 + 12t, 10 + 12t, 7 + 12t, 1)$
    \end{enumerate}
    is transmission irregular.
\end{theorem}

Finally, we give the next seven infinite families of H-trees whose members also represent TI chemical trees of even order.

\begin{theorem}\label{h_tree_th}
    For any $t \in \mathbb{N}_0$, each of the seven H-trees
    \begin{enumerate}[label=\textbf{(\roman*)}]
        \item $\mathrm{H}(1 + t; 6 + 2t, 5 + 2t; 2 + t, 1)$,
        \item $\mathrm{H}(2 + t; 7 + 2t, 6 + 2t; 5 + t, 1)$,
        \item $\mathrm{H}(1 + 2t; 6 + 4t, 5 + 4t; 1 + 2t, 2)$,
        \item $\mathrm{H}(2 + 2t; 8 + 4t, 7 + 4t; 3 + 2t, 1)$,
        \item $\mathrm{H}(3 + 2t; 8 + 4t, 7 + 4t; 6 + 2t, 1)$,
        \item $\mathrm{H}(1 + t; 7 + 4t, 6 + 4t; 6 + 3t, 1)$,
        \item $\mathrm{H}(2; 7 + 4t, 6 + 4t; 5 + 4t, 1)$
    \end{enumerate}
    is transmission irregular.
\end{theorem}

The remainder of the paper shall be structured as follows. In Section \ref{sc_prelim} we will preview some basic theoretical facts mostly regarding the Diophantine equations that will be of use to us later on. Sections \ref{starlike_sc} and \ref{dstarlike_sc} will subsequently serve to provide the proofs of Theorems~\ref{main_th_1} and \ref{main_th_2}, respectively. In the end, Section \ref{sc_families} will be used to disclose the computer-assisted proof of Theorems \ref{x_tree_th} and \ref{h_tree_th} whose Python code is to be found in Appendix \ref{python_xh}.

\section{Preliminaries}\label{sc_prelim}

To begin, we disclose the following folklore auxiliary lemma.

\begin{lemma}\label{neighbor_lemma}
    For any connected graph $G$ and any two of its adjacent vertices $x, y \in V(G)$, the difference $\mathrm{Tr}_G(y) - \mathrm{Tr}_G(x)$ equals the number of vertices that are closer to $x$ than to $y$ minus the number of vertices that are closer to $y$ than to $x$.
\end{lemma}
\begin{proof}
    It suffices to observe that
    \[
        \mathrm{Tr}_G(y) - \mathrm{Tr}_G(x) = \sum_{v \in V(G)} \left( d_G(v, y) - d_G(v, x) \right)
    \]
    and that \[
        d_G(v, y) - d_G(v, x) = \begin{cases}
            1, & \mbox{$v$ is closer to $x$ than to $y$},\\
            -1, & \mbox{$v$ is closer to $y$ than to $x$},\\
            0, & \mbox{otherwise} .
        \end{cases}
    \]
\end{proof}

Now, for each $x \in \mathbb{N}$, let $[x]$ denote the set $\{1, 2, \ldots, x \}$. We state and prove the following short lemma regarding the images of bivariate discrete mappings.

\begin{lemma}\label{trans_lemma}
    For any $C_1, C_2 \in \mathbb{N}$ and $C_3, C_4 \in \mathbb{Z}$, the image of the bivariate mapping $f \colon [C_1] \times [C_2] \to \mathbb{Z}^2$ defined via
    \[
        f(x, y) = (x - y + C_3, x + y + C_4) \qquad \left( x \in [C_1], y \in [C_2] \right)
    \]
    consists of the pairs $(u, v) \in \mathbb{Z}^2$ satisfying the conditions:
    \begin{enumerate}[label=\textbf{(\roman*)}]
        \item $C_3 + C_4 + 2 \le u + v \le C_3 + C_4 + 2 C_1$;
        \item $C_4 - C_3 + 2 \le v - u \le C_4 - C_3 + 2 C_2$;
        \item $u + v \equiv_2 C_3 + C_4$.
    \end{enumerate}
\end{lemma}
\begin{proof}
    Let $W = f([C_1] \times [C_2])$. If we select a fixed $(u, v) \in W$, we then have that $u = x - y + C_3$ and $v = x + y + C_4$ simultaneously hold for some $x \in [C_1]$ and $y \in [C_2]$. From here, it follows that $u + v = C_3 + C_4 + 2x$ and $v - u = C_4 - C_3 + 2y$ and it becomes easy to conclude that all the stated conditions certainly hold. In order to finalize the proof, we have to show that any pair $(u, v) \in \mathbb{Z}^2$ which satisfies the given conditions necessarily belongs to $W$. However, to achieve this, it is enough to notice that for any such $(u, v)$, we have
    \[
        f\left( \frac{u + v - C_3 - C_4}{2}, \frac{v - u + C_3 - C_4}{2} \right) = (u, v) ,
    \]
    where the values $\dfrac{u + v - C_3 - C_4}{2}$ and $\dfrac{v - u + C_3 - C_4}{2}$ are obviously both integers from $[C_1]$ and $[C_2]$, respectively.
\end{proof}

We will now apply Lemma \ref{trans_lemma} in order to obtain an auxiliary lemma on Diophantine equations that will play a key role in proving Theorems \ref{main_th_1} and \ref{main_th_2}.

\begin{lemma}\label{dio_lemma}
    Let $C_1, C_2 \in \mathbb{N}_0$ be such that $C_1 \equiv_2 C_2$. Also, let $C_3 \in \mathbb{Z}$ and $C_4, C_5 \in \mathbb{N}$. The Diophantine equation
    \begin{equation}\label{dio_eq}
        (x^2 + C_1 x) - (y^2 + C_2 y) = C_3
    \end{equation}
    in $(x, y) \in [C_4] \times [C_5]$ has a solution if and only if the integer
    \[
        C_* = C_3 + \frac{(C_1 - C_2)(C_1 + C_2)}{4}
    \]
    has a positive divisor $p \in \mathbb{N}, \, p > \sqrt{\left| C_* \right|}$ such that the following conditions hold:
    \begin{enumerate}[label=\textbf{(\roman*)}]
        \item $C_1 + 2 \le p + \frac{C_*}{p} \le C_1 + 2C_4$;
        \item $C_2 + 2 \le p - \frac{C_*}{p} \le C_2 + 2C_5$;
        \item $p + \frac{C_*}{p} \equiv_2 C_1$.
    \end{enumerate}
\end{lemma}
\begin{proof}
    Given the fact that
    \begin{align*}
        (x^2 + C_1 x) - (y^2 + C_2 y) &= \left( x + \frac{C_1}{2} \right)^2 - \left( y + \frac{C_2}{2} \right)^2 - \frac{C_1^2 - C_2^2}{4}\\
        &= \left( x - y + \frac{C_1 - C_2}{2} \right) \left( x + y + \frac{C_1+C_2}{2} \right) - \frac{C_1^2 - C_2^2}{4} ,
    \end{align*}
    it immediately follows that Eq.\ \eqref{dio_eq} is equivalent to
    \begin{equation}\label{dio_eq_2}
        \left( x - y + \frac{C_1 - C_2}{2} \right) \left( x + y + \frac{C_1+C_2}{2} \right) = C_* .
    \end{equation}
    Now, if we put $q = x - y + \frac{C_1 - C_2}{2} $ and $p = x + y + \frac{C_1 + C_2}{2}$, and then apply Lemma~\ref{trans_lemma}, it becomes evident that the Diophantine equation \eqref{dio_eq_2} has a solution in $(x, y) \in [C_4] \times [C_5]$ if and only if the Diophantine equation
    \begin{equation}\label{dio_eq_3}
        qp = C_*
    \end{equation}
    has a solution $(q, p) \in \mathbb{Z}^2$ such that
    \begin{enumerate}[label=\textbf{(\roman*)}]
        \item $C_1 + 2 \le p + q \le C_1 + 2C_4$;
        \item $C_2 + 2 \le p - q \le C_2 + 2C_5$;
        \item $p + q \equiv_2 C_1$.
    \end{enumerate}
    
    For each solution $(q, p)$ to Eq.\ \eqref{dio_eq_3} where $p \neq 0$, it is obvious that $q = \frac{C_*}{p}$. For this reason, in order to complete the proof, it is enough to show that $p > \sqrt{|C_*|}$ must hold for each solution $(q, p) \in \mathbb{Z}^2$ that satisfies the three aforementioned conditions \textbf{(i)}, \textbf{(ii)} and \textbf{(iii)}. First of all, it is clear that $p > 0$, since if we sum up the inequalities from conditions \textbf{(i)} and \textbf{(ii)}, we directly obtain
    $2p \ge C_1 + C_2 + 2$. On the other hand, if we suppose that $0 < p \le \sqrt{|C_*|}$, then it promptly follows that
    \[
        p - \frac{|C_*|}{p} \le \sqrt{|C_*|} - \frac{|C_*|}{\sqrt{|C_*|}} = 0 ,
    \]
    implying that at least one of the numbers $p + \frac{C_*}{p}, p - \frac{C_*}{p}$ must be nonpositive. However, this is not in accordance with the inequalities given in conditions \textbf{(i)} and \textbf{(ii)}. Thus, $p > \sqrt{|C_*|}$ is necessarily true.
\end{proof}

We end the section by stating another auxiliary lemma to be heavily used later on in Sections~\ref{starlike_sc} and \ref{dstarlike_sc}.

\begin{lemma}\label{func_lemma}
    Let $C_1 \in \mathbb{R}$ and $C_2 > 0$. For any real number $x > \sqrt{|C_1|}$, it is true that
    \[
        x + \frac{C_1}{x} \ge C_2, \qquad \mbox{if and only if} \qquad x \ge g(C_1, C_2),
    \]
    as well as
    \[
        x + \frac{C_1}{x} \le C_2, \qquad \mbox{if and only if} \qquad x \le g(C_1, C_2).
    \]
\end{lemma}
\begin{proof}
    It is straightforwad to see that for any $x > \sqrt{|C_1|}$, the inequality
    \[
        x + \frac{C_1}{x} \ge C_2
    \]
    is equivalent to
    \[
        x^2 - C_2 x + C_1 \ge 0,
    \]
    which is further equivalent to
    \begin{equation}\label{paux_1}
        \left( x - \frac{C_2}{2} \right)^2 \ge \frac{C_2^2 - 4C_1}{4} .
    \end{equation}
    Now, if $C_2^2 - 4C_1 < 0$, then $g(C_1, C_2) = -\infty$ and it follows that Eq.\ \eqref{paux_1} and $x \ge g(C_1, C_2)$ both necessarily hold. On the other hand, if $C_2^2 - 4C_1 \ge 0$, then Eq.~\eqref{paux_1} directly becomes equivalent to
    \begin{equation}\label{paux_2}
        x \le \zeta \quad \mbox{or} \quad x \ge g(C_1, C_2) ,
    \end{equation}
    where $\zeta = \frac{1}{2} \left( C_2 - \sqrt{C_2^2 - 4 C_1} \right)$.
    However, due to $\zeta \, g(C_1, C_2) = C_1$, it follows that either $C_1 \le 0$, in which case $\zeta \le 0$, or $C_1 > 0$, thus implying $g(C_1, C_2) \ge \zeta > 0$ and
    \begin{align*}
        \zeta^2 \le \zeta \, g(C_1, C_2) = C_1 .
    \end{align*}
    Either way, $\zeta \le \sqrt{|C_1|}$ must be true, which means that Eq.\ \eqref{paux_2} is surely equivalent to just $x \ge g(C_1, C_2)$.

    Similarly, for any $x > \sqrt{|C_1|}$, the inequality
    \[
        x + \frac{C_1}{x} \le C_2
    \]
    is equivalent to
    \[
        x^2 - C_2 x + C_1 \le 0,
    \]
    which is then equivalent to
    \begin{equation}\label{paux_3}
        \left( x - \frac{C_2}{2} \right)^2 \le \frac{C_2^2 - 4C_1}{4} .
    \end{equation}
    If $C_2^2 - 4C_1 < 0$, then $g(C_1, C_2) = -\infty$ and it is clear that Eq.\ \eqref{paux_3} and $x \le g(C_1, C_2)$ are surely both not satisfied. Finally, if we have $C_2^2 - 4C_1 \ge 0$, then Eq.\ \eqref{paux_3} quickly becomes equivalent to
    \begin{equation}\label{paux_4}
        \zeta \le x \le g(C_1, C_2) ,
    \end{equation}
    where $\zeta = \frac{1}{2} \left( C_2 - \sqrt{C_2^2 - 4 C_1} \right)$. We have already shown that $\zeta \le \sqrt{|C_1|}$ certainly holds, which immediately implies that Eq.\ \eqref{paux_4} is equivalent to $x \le g(C_1, C_2)$.
\end{proof}

\section{Starlike trees}\label{starlike_sc}

In this section, we will provide the proof of Theorem \ref{main_th_1}. Let $T = \mathrm{S}(A_1, A_2, \ldots, A_k)$ be a given starlike tree whose vertex of degree $k$ is denoted by $v_{10} \equiv v_{20} \equiv \cdots \equiv v_{k0}$ and whose remaining vertices we shall signify by $v_{ij}$, $1 \le i \le k, \, 1 \le j \le A_i$, so that the path $v_{i0} v_{i1} v_{i2} \cdots x_{iA_i}$ represents a pendent path of length $A_i$. For starters, if we have $A_i = A_j$ for some $1 \le i < j \le k$, it is then straightforward to observe that there necessarily exists an automorphism of $T$ which maps $v_{i1}$ to $v_{j1}$. Hence, these two vertices must have the same transmission and the graph $T$ cannot be TI. This leads us to the following proposition.

\begin{proposition}\label{prop_1}
    If the values $A_1, A_2, \ldots, A_k$ are not all mutually distinct, then the starlike tree $\mathrm{S}(A_1, A_2, \ldots, A_k)$ is not transmission irregular.
\end{proposition}

Furthermore, from Lemma \ref{neighbor_lemma} we promptly get
\begin{equation}\label{aux_1}
    \mathrm{Tr}(v_{ij}) - \mathrm{Tr}(v_{i, j-1}) = n - 2A_i + 2j - 2 \qquad (1 \le i \le k, 1 \le j \le A_i) ,
\end{equation}
which immediately leads us to 
\[
    \mathrm{Tr}(v_{ij}) = \mathrm{Tr}(v_{i0}) + \sum_{h=1}^{j} (n - 2A_i + 2h - 2) .
\]
Therefore, for any $1 \le i \le k, \, 0 \le j \le A_i$, we have
\begin{align}\nonumber
    \mathrm{Tr}(v_{ij}) &= \mathrm{Tr}(v_{i0}) + j \, \frac{(n - 2A_i) + (n - 2A_i + 2j - 2)}{2}\\
    \label{aux_4} &= \mathrm{Tr}(v_{i0}) + j(n - 2A_i + j - 1) .
\end{align}
By relying on Eq.\ \eqref{aux_4}, we are now able to prove the following proposition.

\begin{proposition}\label{prop_2}
    If $A_i \ge \frac{n}{2}$ holds for some $1 \le i \le k$, then the starlike tree $\mathrm{S}(A_1, A_2, \ldots, \linebreak A_k)$ is not transmission irregular.
\end{proposition}
\begin{proof}
    If $A_i \ge \frac{n}{2}$, we have $2A_i + 1 - n \ge 1$ together with $2A_i + 1 - n \le A_i$. Thus, Eq.~\eqref{aux_4} allows us to compute $\mathrm{Tr}(v_{i, 2 A_i + 1 - n}) = \mathrm{Tr}(v_{i0})$, hence the starlike tree is certainly not TI.
\end{proof}

By virtue of Propositions \ref{prop_1} and \ref{prop_2}, we see that in case that at least one of the conditions \textbf{(i)} or \textbf{(ii)} from Theorem \ref{main_th_1} does not hold, then the given starlike tree fails to be TI, as desired. In order to prove the validity of the theorem, we will suppose that these two conditions indeed are satisfied, then show that $\mathrm{S}(A_1, A_2, \ldots, A_k)$ is TI if and only if condition \textbf{(iii)} holds.

Provided that conditions \textbf{(i)} and \textbf{(ii)} do hold, Eq.\ \eqref{aux_1} dictates that the finite sequence
\[
    \mathrm{Tr}(v_{i0}), \mathrm{Tr}(v_{i1}), \mathrm{Tr}(v_{i2}), \ldots, \mathrm{Tr}(v_{iA_i})
\]
must be strictly increasing for each $1 \le i \le k$. For this reason, the given starlike tree is not TI if and only if there exist some $1 \le i < j \le k$ and $1 \le x \le A_i, 1 \le y \le A_j$ such that $\mathrm{Tr}(v_{ix}) = \mathrm{Tr}(v_{jy})$. However, Eq.\ \eqref{aux_4} implies that inspecting the given equality gets down to solving the Diophantine equation
\begin{equation}\label{aux_5}
    x(n - 2A_i + x - 1) = y(n - 2A_j + y - 1)
\end{equation}
in $(x, y) \in [A_i] \times [A_j]$. Given the fact that Eq.\ \eqref{aux_5} is clearly equivalent to
\[
    (x^2 + x(n - 2A_i - 1)) - (y^2 + y(n - 2A_j - 1)) = 0 ,
\]
where $n - 2A_i - 1 \ge 0$ and $n - 2A_j - 1 \ge 0$, Lemma \ref{dio_lemma} tells us that Eq.\ \eqref{aux_5} has a solution if and only if the integer
\begin{align*}
    &\frac{((n - 2A_i - 1) - (n - 2A_j - 1))((n - 2A_i - 1) + (n - 2A_j - 1))}{4} =\\
    &\qquad\qquad = \frac{(2A_j - 2A_i)(2n - 2 - 2A_i - 2A_j)}{4} = \alpha_{ij}
\end{align*}
has a positive divisor $p \in \mathbb{N}, \, p > \sqrt{|\alpha_{ij}|}$ such that the following conditions hold:
\begin{enumerate}[label=\textbf{(\roman*)}]
    \item $n + 1 - 2A_i \le p + \frac{\alpha_{ij}}{p} \le n - 1$;
    \item $n + 1 - 2A_j \le p - \frac{\alpha_{ij}}{p} \le n - 1$;
    \item $p + \frac{\alpha_{ij}}{p} \equiv_2 n - 1$.
\end{enumerate}
This quickly brings us to the following proposition.

\begin{proposition}\label{prop_3}
    If conditions \textbf{(i)} and \textbf{(ii)} from Theorem \ref{main_th_1} are satisfied, then the starlike tree $\mathrm{S}(A_1, A_2, \ldots, A_k)$ is TI if and only if there do not exist any $1 \le i < j \le k$ for which there is a positive divisor $p \in \mathbb{N}, \, p > \sqrt{|\alpha_{ij}|}$ of $\alpha_{ij}$ that satisfies the conditions:
    \begin{enumerate}[label=\textbf{(\roman*)}]
        \item $n + 1 - 2A_i \le p + \frac{\alpha_{ij}}{p} \le n - 1$;
        \item $n + 1 - 2A_j \le p - \frac{\alpha_{ij}}{p} \le n - 1$;
        \item $p + \frac{\alpha_{ij}}{p} \equiv_2 n - 1$.
    \end{enumerate}
\end{proposition}
Condition \textbf{(iii)} from Proposition \ref{prop_3} clearly coincides with condition \textbf{(a)} from Theorem~\ref{main_th_1}. For this reason, in order to finalize the proof, it is enough to show that the other two conditions \textbf{(i)} and \textbf{(ii)} are equivalent to conditions \textbf{(b)} and \textbf{(c)} from Theorem~\ref{main_th_1}. This is precisely what we will demonstrate in the remainder of the section, thereby completing the proof as follows.

\bigskip\noindent
\emph{Proof of Theorem \ref{main_th_1}}.\quad
If at least one of the conditions \textbf{(i)} or \textbf{(ii)} fails to hold, then Propositions \ref{prop_1} and \ref{prop_2} guarantee that the starlike tree $\mathrm{S}(A_1, A_2, \ldots, A_k)$ is not TI, as desired. Suppose that both of these conditions are satisfied. In this case, Proposition \ref{prop_3} provides the precise conditions under which the given tree is TI. We will now demonstrate that these conditions are equivalent to those given in the theorem.

It is obvious that condition \textbf{(iii)} from Proposition \ref{prop_3} is equivalent to condition~\textbf{(a)} from Theorem \ref{main_th_1}. Furthermore, Lemma \ref{func_lemma} guarantees that condition~\textbf{(i)} from Proposition~\ref{prop_3} is equivalent to condition~\textbf{(b)} from Theorem~\ref{main_th_1}. Finally, bearing in mind that $-\alpha_{ij} = \alpha_{ji}$, Lemma \ref{func_lemma} also dictates that condition \textbf{(ii)} from Proposition \ref{prop_3} must be equivalent to condition \textbf{(c)} from Theorem \ref{main_th_1}. \hfill\qed

\bigskip

\begin{remark}
    We also point out that, provided conditions \textbf{(i)} and \textbf{(ii)} from Theorem~\ref{main_th_1} hold, none of the values
    \[
        g(\alpha_{ij}, n + 1 - 2A_i), g(\alpha_{ij}, n-1), g(\alpha_{ji}, n + 1 - 2A_j), g(\alpha_{ji}, n-1)
    \]
    from condition \textbf{(iii)} will actually be equal to $-\infty$. In order to verify this, it is sufficient to consider the function $f \colon (0, +\infty) \to \mathbb{R}$ defined via
    \[
        f(x) = x + \frac{\alpha_{ij}}{x}
    \]
    for some $1 \le i, j \le k, \, i \neq j$. From
    \begin{align*}
        (n - 1 - A_i - A_j) - (A_j - A_i) = n - 1 - 2A_j &\ge 0,\\
        (n - 1 - A_i - A_j) - (A_i - A_j) = n - 1 - 2A_i &\ge 0,
    \end{align*}
    we immediately get $n - 1 - A_i - A_j \ge |A_j - A_i | > 0$, which leads us to
    \[
        f(n-1-A_i-A_j) = (n-1-A_i-A_j) + (A_j - A_i) = n - 1 - 2A_i .
    \]
    Since $n - 1 - 2A_i < n + 1 - 2A_i \le n - 1$, from the continuity of $f$ and its limit property $\displaystyle\lim_{x \to +\infty}f(x) = +\infty$, we conclude that $n + 1 - 2A_i$ and $n-1$ both belong to $f((0, +\infty))$. From here, it follows that the corresponding quadratic equations
    \begin{align*}
        x^2 - x(n + 1 - 2A_i) + \alpha_{ij} &= 0,\\
        x^2 - x(n - 1) + \alpha_{ij} &= 0,
    \end{align*}
    necessarily have real solutions. \hfill\qed
\end{remark}

\section{Double starlike trees}\label{dstarlike_sc}

Here, we will give the proof of Theorem \ref{main_th_2}. Let $T = \mathrm{DS}(C; A_1, A_2, \ldots, A_k; B_1, B_2, \ldots, B_m)$ be a given double starlike tree whose vertices we shall denote as follows:
\begin{enumerate}[label=\textbf{(\roman*)}]
    \item the vertex of degree $k$ which represents the starting point of the pendent paths of lengths $A_1, \ldots, A_k$ will be denoted by $v_{00} \equiv v_{10} \equiv v_{20} \equiv \cdots \equiv v_{k0}$;
    \item the $k$ pendent paths starting at $v_{00}$ will bear the form $v_{i0} v_{i1} v_{i2} \cdots v_{iA_i}$, for $1 \le i \le k$;
    \item the other vertex of degree greater than two that is different from $v_{00}$ will be denoted by $v_{0C} \equiv u_{10} \equiv u_{20} \equiv \cdots \equiv u_{m0}$;
    \item the path from $v_{00}$ to $v_{0C}$ will be written as $v_{00} v_{01} v_{02} \cdots v_{0C}$;
    \item the $m$ pendent paths starting at $v_{0C}$ will be $u_{i0} u_{i1} u_{i2} \cdots u_{iB_i}$, for $1 \le i \le m$.
\end{enumerate}

To begin, if we have that $A_i = A_j$ holds for some $1 \le i < j \le k$, it is trivial to notice that there must exist an automorphism that maps $v_{i1}$ to $v_{j1}$. From here, it promptly follows that these two vertices surely have the same transmission and that the tree $T$ cannot be TI. An analogous conclusion can be made regarding the values $B_1, B_2, \ldots, B_m$, which brings us to the next proposition.

\begin{proposition}\label{dprop_1}
    If the values $A_1, A_2, \ldots, A_k$ are not all mutually distinct or if the values $B_1, B_2, \ldots, B_m$ are not all mutually distinct, then the double starlike tree $\mathrm{DS}(C; A_1, A_2, \linebreak \ldots, A_k; B_1, B_2, \ldots, B_m)$ is not transmission irregular.
\end{proposition}

With the help of Lemma \ref{neighbor_lemma}, it is easy to see that
\begin{equation}\label{aux_7}
    \mathrm{Tr}(v_{ij}) - \mathrm{Tr}(v_{i, j-1}) = n - 2A_i + 2j - 2 \qquad (1 \le i \le k, 1 \le j \le A_i) .
\end{equation}
By using the same strategy as in Section \ref{starlike_sc}, we are now able to obtain
\begin{equation}\label{aux_2}
    \mathrm{Tr}(v_{ij}) = \mathrm{Tr}(v_{i0}) + j(n - 2A_i + j - 1) \qquad (1 \le i \le k, 0 \le j \le A_i) .
\end{equation}
Moreover, it can be analogously shown that
\begin{equation}\label{aux_3}
    \mathrm{Tr}(u_{ij}) = \mathrm{Tr}(u_{i0}) + j(n - 2B_i + j - 1) \qquad (1 \le i \le m, 0 \le j \le B_i) .
\end{equation}
Another implementation of Lemma \ref{neighbor_lemma} leads to
\begin{equation}\label{aux_8}
    \mathrm{Tr}(v_{0j}) - \mathrm{Tr}(v_{0, j-1}) = 2(1 + A_*) - n + 2j - 2 \qquad (1 \le j \le C) ,
\end{equation}
which then immediately implies
\[
    \mathrm{Tr}(v_{0j}) = \mathrm{Tr}(v_{00}) + \sum_{h=1}^j \left( 2(1 + A_*) - n + 2h - 2 \right)
\]
for any $0 \le j \le C$. Therefore,
\begin{align}
    \nonumber \mathrm{Tr}(v_{0j}) &= \mathrm{Tr}(v_{00}) + j \, \frac{\left( 2(1 + A_*) - n \right) + \left( 2(1 + A_*) - n + 2j - 2\right)}{2}\\
    \label{aux_6} &= \mathrm{Tr}(v_{00}) + j \left( 2(1 + A_*) - n + j - 1 \right) .
\end{align}
By using Eqs.\ \eqref{aux_2}, \eqref{aux_3} and \eqref{aux_6}, we are able to concisely demonstrate the validity of the next proposition.

\begin{proposition}\label{dprop_2}
    If $A_i \ge \frac{n}{2}$ holds for some $1 \le i \le k$, or if $B_i \ge \frac{n}{2}$ holds for some $1 \le i \le m$, or if $1 + A_* \le \frac{n}{2}$ is true, then the double starlike tree $\mathrm{DS}(C; A_1, A_2, \ldots, A_k; B_1, B_2, \ldots, \linebreak B_m)$ is not transmission irregular.
\end{proposition}
\begin{proof}
    If $A_i \ge \frac{n}{2}$, it is straightforward to see that $2A_i + 1 - n \ge 1$ and $2A_i + 1 - n \le A_i$, hence it immediately follows that $\mathrm{Tr}(v_{i, 2A_i + 1 - n}) = \mathrm{Tr}(v_{i0})$, by virtue of Eq.\ \eqref{aux_2}. Thus, the tree is not TI. An analogous logical reasoning can be used to apply Eq.\ \eqref{aux_3} to resolve the case when $B_i \ge \frac{n}{2}$ for some $1 \le i \le m$. On the other hand, if we suppose that $1 + A_* \le \frac{n}{2}$ is satisfied, we then have $n + 1 - 2(1+A_*) \ge 1$, as well as
    \[
        n + 1 - 2(1 + A_*) = \left(n - 1 - A_* \right) - A_* = C + B_* - A_* \le C .
    \]
    By using Eq.\ \eqref{aux_6}, we directly obtain $\mathrm{Tr}(v_{0, n + 1 - 2(1+A_*)}) = \mathrm{Tr}(v_{00})$, which means that the given double starlike tree cannot be TI.
\end{proof}

Propositions \ref{dprop_1} and \ref{dprop_2} tell us that if at least one of the conditions \textbf{(i)} or \textbf{(ii)} from Theorem \ref{main_th_2} is not satisfied, then the given double starlike tree cannot be TI. For this reason, in order to complete the proof of Theorem~\ref{main_th_2}, it is enough to suppose that conditions \textbf{(i)} and \textbf{(ii)} hold, then demonstrate that $\mathrm{DS}(C; A_1, A_2, \ldots, A_k; B_1, B_2, \ldots, B_m)$ is TI if and only if conditions \textbf{(iii)}, \textbf{(iv)}, \textbf{(v)} and \textbf{(vi)} are all satisfied.

In the case that conditions \textbf{(i)} and \textbf{(ii)} from Theorem \ref{main_th_2} hold, Eq.\ \eqref{aux_7} quickly implies that the finite sequence
\[
    \mathrm{Tr}(v_{i0}), \mathrm{Tr}(v_{i1}), \mathrm{Tr}(v_{i2}), \ldots, \mathrm{Tr}(v_{iA_i})
\]
must be strictly increasing for each $1 \le i \le k$, and we may analogously conclude that the finite sequence
\[
    \mathrm{Tr}(u_{i0}), \mathrm{Tr}(u_{i1}), \mathrm{Tr}(u_{i2}), \ldots, \mathrm{Tr}(u_{iB_i})
\]
is necessarily strictly increasing for each $1 \le i \le m$. Furthermore, Eq.\ \eqref{aux_8} tells us that the finite sequence
\[
    \mathrm{Tr}(v_{00}), \mathrm{Tr}(v_{01}), \mathrm{Tr}(v_{02}), \ldots, \mathrm{Tr}(v_{0C})
\]
is also certainly strictly increasing. Bearing all of this in mind, it is not difficult to observe that the given double starlike tree is not TI if and only if at least one equation in $(x, y) \in \mathbb{Z}^2$ is solvable from the following four sets:
\begin{alignat}{3}
    \label{ugly_1} \mathrm{Tr}(v_{ix}) &= \mathrm{Tr}(v_{jy}) \qquad && (1 \le x \le A_i, 1 \le y \le A_j) \qquad && (1 \le i < j \le k),\\
    \label{ugly_2} \mathrm{Tr}(u_{ix}) &= \mathrm{Tr}(u_{jy}) \qquad && (1 \le x \le B_i, 1 \le y \le B_j) \qquad && (1 \le i < j \le m),\\
    \label{ugly_3} \mathrm{Tr}(v_{ix}) &= \mathrm{Tr}(v_{0y}) \qquad && (1 \le x \le A_i, 1 \le y \le C) \qquad && (1 \le i \le k),\\
    \label{ugly_4} \mathrm{Tr}(v_{ix}) &= \mathrm{Tr}(u_{jy}) \qquad && (1 \le x \le A_i, 1 \le y \le B_j) \qquad && (1 \le i \le k, 1 \le j \le m) .
\end{alignat}
Now, Eq.\ \eqref{aux_6} promptly leads to
\begin{align*}
    \mathrm{Tr}(v_{0C}) &= \mathrm{Tr}(v_{00}) + C(2(1 + A_*) - n + C - 1)\\
    &= \mathrm{Tr}(v_{00}) + C((1 + A_*) - (n - C - A_* ))\\
    &= \mathrm{Tr}(v_{00}) + C((1 + A_*) - (1 + B_*))\\
    &= \mathrm{Tr}(v_{00}) + C(A_* - B_*) ,
\end{align*}
which then together with Eqs.\ \eqref{aux_2}, \eqref{aux_3} and \eqref{aux_6} helps us conclude that Eqs.\ \eqref{ugly_1}, \eqref{ugly_2}, \eqref{ugly_3} and \eqref{ugly_4} are equivalent to the next four equations:
\begin{align*}
    (x^2 + (n - 2A_i - 1)x) - (y^2 + (n - 2A_j - 1)y) &= 0,\\
    (x^2 + (n - 2B_i - 1)x) - (y^2 + (n - 2B_j - 1)y) &= 0,\\
    (x^2 + (n - 2A_i - 1)x) - (y^2 + (2A_* - n + 1)y) &= 0,\\
    (x^2 + (n - 2A_i - 1)x) - (y^2 + (n - 2B_j - 1)y) &= C(A_* - B_*),
\end{align*}
respectively. By using Lemma \ref{dio_lemma} in the same way as it was done in Section \ref{starlike_sc}, a quick computation is all it takes to reach the following proposition.

\begin{proposition}\label{dprop_3}
    Provided conditions \textbf{(i)} and \textbf{(ii)} from Theorem \ref{main_th_2} are satisfied, the double starlike tree $\mathrm{DS}(C; A_1, A_2, \ldots, A_k; B_1, B_2, \ldots, B_m)$ is transmission irregular if and only if the following conditions hold:
    \begin{enumerate}[label=\textbf{(\roman*)}]
        \item there do not exist any $1 \le i < j \le k$ for which there is a positive divisor $p \in \mathbb{N}, \, p > \sqrt{|\alpha_{ij}|}$ of $\alpha_{ij}$ that satisfies the conditions:
        \begin{enumerate}[label=\textbf{(\alph*)}]
            \item $n + 1 - 2A_i \le p + \frac{\alpha_{ij}}{p} \le n - 1$;
            \item $n + 1 - 2A_j \le p - \frac{\alpha_{ij}}{p} \le n - 1$;
            \item $p + \frac{\alpha_{ij}}{p} \equiv_2 n - 1$;
        \end{enumerate}
        \item there do not exist any $1 \le i < j \le m$ for which there is a positive divisor $p \in \mathbb{N}, \, p > \sqrt{|\beta_{ij}|}$ of $\beta_{ij}$ that satisfies the conditions:
        \begin{enumerate}[label=\textbf{(\alph*)}]
            \item $n + 1 - 2B_i \le p + \frac{\beta_{ij}}{p} \le n - 1$;
            \item $n + 1 - 2B_j \le p - \frac{\beta_{ij}}{p} \le n - 1$;
            \item $p + \frac{\beta_{ij}}{p} \equiv_2 n - 1$;
        \end{enumerate}
        \item there does not exist any $1 \le i \le k$ for which there is a positive divisor $p \in \mathbb{N}, \, p > \sqrt{|\gamma_{i}|}$ of $\gamma_{i}$ that satisfies the conditions:
        \begin{enumerate}[label=\textbf{(\alph*)}]
            \item $n + 1 - 2A_i \le p + \frac{\gamma_{i}}{p} \le n - 1$;
            \item $2A_* - n + 3 \le p - \frac{\gamma_{i}}{p} \le 2A_* + 2C - n + 1$;
            \item $p + \frac{\gamma_{i}}{p} \equiv_2 n - 1$;
        \end{enumerate}
        \item there do not exist any $1 \le i \le k, 1 \le j \le m$ for which there is a positive divisor $p \in \mathbb{N}, \, p > \sqrt{|\delta_{ij}|}$ of $\delta_{ij}$ that satisfies the conditions:
        \begin{enumerate}[label=\textbf{(\alph*)}]
            \item $n + 1 - 2A_i \le p + \frac{\delta_{ij}}{p} \le n - 1$;
            \item $n + 1 - 2B_j \le p - \frac{\delta_{ij}}{p} \le n - 1$;
            \item $p + \frac{\delta_{ij}}{p} \equiv_2 n - 1$.
        \end{enumerate}
    \end{enumerate}
\end{proposition}

It now becomes straightforward to finalize the proof of Theorem \ref{main_th_2} in the same manner as it was done with Theorem \ref{main_th_1} in Section \ref{starlike_sc}.

\bigskip\noindent
\emph{Proof of Theorem \ref{main_th_2}}.\quad
If at least one of the conditions \textbf{(i)} or \textbf{(ii)} is not safisfied, then Propositions \ref{dprop_1} and \ref{dprop_2} imply that the double starlike tree $\mathrm{DS}(C; A_1, A_2, \ldots, A_k; B_1, B_2, \linebreak \ldots, B_m)$ is not TI. Now, suppose that both of these conditions hold. In this scenario, Proposition \ref{dprop_3} dictates the precise conditions under which the tree of interest is TI. However, a direct application of Lemma \ref{func_lemma} directly shows that these conditions are equivalent to the corresponding ones stated in Theorem \ref{main_th_2}. \hfill\qed

\section{Families of X-trees and H-trees}\label{sc_families}

In the last section we will disclose a mechanism that can be used to prove Theorems~\ref{x_tree_th} and \ref{h_tree_th} by using Theorems \ref{main_th_1} and \ref{main_th_2}, respectively. We will first briefly explain the logic behind the strategy and how it can be applied to each of the 11 families of graphs given in the two theorems. Afterwards, we will end the section by providing the full proof concerning the $\mathrm{X}(7 + 12t, 6 + 12t, 3 + 12t, 1)$ and $\mathrm{H}(1 + t; 6 + 2t, 5 + 2t; 2 + t, 1)$ infinite families. We omit the classic proof regarding the remaining nine families and instead provide the computer-assisted proof that automatically handles all $11$ families and whose Python implementation can be found in Appendix \ref{python_xh}.

To begin, it is straightforward to check that all the graphs appearing in Theorems \ref{x_tree_th} and \ref{h_tree_th} are of even order. Furthermore, it is trivial to verify that conditions \textbf{(i)} and \textbf{(ii)} from Theorem \ref{main_th_1} (resp.\ conditions \textbf{(i)} and \textbf{(ii)} from Theorem \ref{main_th_2}) are both satisfied for each of the aforementioned X-trees (resp.\ H-trees). For this reason, in order to complete the proof, it is sufficient to focus only on condition \textbf{(iii)} from Theorem \ref{main_th_1} (resp.\ conditions \textbf{(iii)}, \textbf{(iv)}, \textbf{(v)} and \textbf{(vi)} from Theorem~\ref{main_th_2}). However, having in mind a fixed family out of the 11 given ones, all of these conditions get down to inspecting whether for any $t_0 \in \mathbb{N}_0$, there exists no $p \in \mathbb{N}, \, p \mid F(t_0), \, p > \sqrt{|F(t_0)|}$ satisfying the three subconditions \textbf{(a)}, \textbf{(b)} and \textbf{(c)}, where $F(t) \in \mathbb{Z}[t], \, \deg F(t) \le 2$ is some according polynomial. Here, it is worth pointing out that subcondition \textbf{(a)} simply dictates that $p$ and $\frac{F(t_0)}{p}$ must be of different parities, since all the graphs of interest are of even order.

It can be manually verified that while applying Theorem \ref{main_th_1} or \ref{main_th_2} on each of the $11$ given families, all the $g(C_1, C_2)$ expressions that appear throughout subconditions \textbf{(b)} and \textbf{(c)} are such that $g(C_1, C_2) \neq -\infty$, i.e., $C_2^2 - 4C_1 \ge 0$. Moreover, we can regard $C_2^2 - 4C_1$ as a $\mathbb{Z}[t]$ polynomial with the degree of at most two and it can be checked that this polynomial has nonnegative coefficients. For this reason, it becomes easy to approximate the value $\sqrt{C_2^2 - 4 C_1}$ via linear polynomials. Thus, the expressions from subconditions \textbf{(b)} and \textbf{(c)} can be used to obtain an inequality bearing the form
\begin{equation}\label{lower_upper_ineq}
    l_0 + l_1 t_0 \le p \le u_0 + u_1 t_0 \qquad (t_0 \in \mathbb{N}_0)
\end{equation}
for some fixed $l_0, l_1, u_0, u_1 \in \mathbb{N}$. If it can be guaranteed that $l_0 + l_1 t_0 > u_0 + u_1 t_0$ holds for any $t_0 \in \mathbb{N}_0$, then there exists no integer $p$ between these two values and there is nothing left to discuss. Otherwise, it becomes convenient to divide the proof into two cases.

\bigskip\noindent
\emph{Case 1: $\deg F(t) \le 1$}.\quad
This case can be resolved within the following five steps.

\begin{enumerate}[label=\textbf{(\arabic*)}]
    \item Verify that $F(t_0) \neq 0$.
    \item If we have that $l_0 + l_1 t_0 > |F(t_0)| > 0$, it is clear that no suitable divisor $p$ can exist. Otherwise, verify that $|F(t_0)| > u_0 + u_1 t_0$ necessarily holds, which immediately implies $\frac{|F(t_0)|}{p} \ge 2$.
    \item If $|F(t_0)| < 2(l_0 + l_1 t_0)$, then $\frac{|F(t_0)|}{p} \ge 2$ cannot hold, hence there exists no desired $p \in \mathbb{N}$. Otherwise, verify that all the $F(t)$ polynomial coefficients are divisible by four. This guarantees that if $\frac{F(t_0)}{p} = 2$, then $p$ is surely even, which implies that subcondition \textbf{(a)} is not satisfied. For this reason, we may further assume that $\frac{F(t_0)}{p} \ge 3$.
    \item If $|F(t_0)| < 3(l_0 + l_1 t_0)$, then $\frac{|F(t_0)|}{p} \ge 3$ cannot be true, which tells us that no required $p$ can exist. Otherwise, verify that $F(t)$ is such that all of its coefficients are divisible by three except for the free term. Hence, $3 \nmid F(t_0)$, which leads us to $\frac{F(t_0)}{p} \ge 4$.
    \item Verify that $|F(t_0)| < 4(l_0 + l_1 t_0)$. This means that $\frac{|F(t_0)|}{p} \ge 4$ certainly does not hold, which guarantees that no suitable $p$ exists.
\end{enumerate}

\bigskip\noindent
\emph{Case 2: $\deg F(t) = 2$}.\quad
In this scenario, the first step is to verify that Eq.~\eqref{lower_upper_ineq} can be obtained in such a way that $l_1 = u_1$. This leads us to
\[
    p \in \left\{ l_1 t_0 + l_0, l_1 t_0 + (l_0 + 1), \ldots, l_1 t_0 + (u_0 - 1), l_1 t_0 + u_0 \right\} ,
\]
which guarantees that $p \in \mathbb{N}$ could have only finitely many potential values. Each of these options should be inspected and discarded separately. It is also important to check that the leading coefficient of $F(t)$ is necessarily divisible by $l_1$. This observation makes it convenient to discard each subcase $p = l_1 t_0 + \theta$, where $\theta \in \mathbb{N}, \, l_0 \le \theta \le u_0$, by applying the polynomial division on $F(t)$ by $l_1 t + \theta$. The remainder becomes straightforward to handle and by performing some relatively simple computations afterwards, it is not difficult to show that either $p \nmid F(t_0)$, or $p \mid F(t_0)$ and $p \equiv_2 \frac{F(t_0)}{p}$. In either scenario, we have that no suitable positive divisor $p$ exists.

\bigskip
In the remainder of the section, we will implement the aforementioned strategy in order to prove that $\mathrm{X}(7 + 12t, 6 + 12t, 3 + 12t, 1)$ and $\mathrm{H}(1 + t; 6 + 2t, 5 + 2t; 2 + t, 1)$ are both TI for any $t \in \mathbb{N}_0$.

\bigskip\noindent
\emph{Proof of Theorem \ref{x_tree_th}, \textbf{(i)}}.\quad
Let $t \in \mathbb{N}_0$ be fixed. We clearly have
\begin{alignat*}{4}
    A_1 &= 7 + 12t, & \qquad A_2 &= 6 + 12t, & \qquad A_3 &= 3 + 12t, & \qquad A_4 &= 1,\\
    & & \qquad n &= 18 + 36t, & \qquad \frac{n}{2} &= 9 + 18t, &
\end{alignat*}
hence it is obvious that the graph order is even and that conditions \textbf{(i)} and \textbf{(ii)} from Theorem \ref{main_th_1} are surely satisfied. In order to complete the proof, it is sufficient to show that condition \textbf{(iii)} also holds by splitting the problem into six corresponding cases.

\bigskip\noindent
\emph{Case 1: $p \mid \alpha_{12}$}.\quad
In this case, we have
\begin{align*}
    \alpha_{12} &= \left( (18 + 36t) - 1 - (7 + 12t) - (6 + 12t)\right)\left( (6 + 12t) - (7 + 12t)\right)\\
    &= (4 + 12t)(-1) = -12t - 4,
\end{align*}
together with $\alpha_{21} = 12t + 4$. Bearing in mind that
\begin{align*}
    g(\alpha_{21}, n + 1 - 2A_2) &= \frac{1}{2} \left( 12t + 7 + \sqrt{(12t+7)^2 - 4(12t+4)}\right)\\
    &= \frac{1}{2} \left( 12t + 7 + \sqrt{144t^2 + 120t + 33}\right)\\
    &= \frac{1}{2} \left( 12t + 7 + \sqrt{(12t+5)^2 + 8}\right)\\
    &> \frac{1}{2} \left( 12t + 7 + 12t + 5 \right) = 12t + 6,
\end{align*}
it becomes evident that if there were a $p \in \mathbb{N}$ satisfying subcondition \textbf{(c)}, then we would have $p > |\alpha_{12}| > 0$, hence such a $p$ could not be a divisor of $\alpha_{12}$.

\bigskip\noindent
\emph{Case 2: $p \mid \alpha_{13}$}.\quad
Here, it is easy to see that
\begin{align*}
    \alpha_{13} &= \left( (18 + 36t) - 1 - (7 + 12t) - (3 + 12t)\right)\left( (3 + 12t) - (7 + 12t)\right)\\
    &= (7 + 12t)(-4) = -48t - 28 ,
\end{align*}
as well as $\alpha_{31} = 48t + 28$. Since
\begin{align*}
    g(\alpha_{31}, n+1-2A_3) &= \frac{1}{2} \left( 12t + 13 + \sqrt{(12t+13)^2 - 4(48t+28)}\right)\\
    &= \frac{1}{2} \left( 12t + 13 + \sqrt{144t^2 + 120t + 57}\right)\\
    &= \frac{1}{2} \left( 12t + 13 + \sqrt{(12t+5)^2 + 32}\right)\\
    &> \frac{1}{2} \left( 12t + 13 + 12t + 5 \right) = 12t + 9
\end{align*}
and
\begin{align*}
    g(\alpha_{31}, n-1) &= \frac{1}{2} \left( 36t + 17 + \sqrt{(36t+17)^2 - 4(48t+28)}\right)\\
    &= \frac{1}{2} \left( 36t + 17 + \sqrt{1296t^2 + 1032t + 177}\right)\\
    &= \frac{1}{2} \left( 36t + 17 + \sqrt{(36t + 15)^2 - 48t - 48}\right)\\
    &< \frac{1}{2} \left( 36t + 17 + 36t + 15 \right) = 36t + 16 ,
\end{align*}
if we suppose that there exists a $p \in \mathbb{N}, \, p \mid \alpha_{13}, \, p > \sqrt{|\alpha_{13}|}$ such that the subconditions \textbf{(a)}, \textbf{(b)} and \textbf{(c)} all hold, we immediately obtain
\begin{equation}\label{xaux_1}
    12t + 10 \le p \le 36t + 15 .
\end{equation}

If $\frac{\alpha_{31}}{p} = 1$, then $p = 48t + 28$, which contradicts Eq.~\eqref{xaux_1}. On the other hand, if $\frac{\alpha_{31}}{p} = 2$, we then get $p = 24t + 14$, which leads to a contradiction due to subcondition~\textbf{(a)}. It is trivial to check that $3 \nmid \alpha_{31}$, which further implies $\frac{\alpha_{31}}{p} \ge 4$. However, we now obtain $p \le 12t + 7$, which contradicts Eq.~\eqref{xaux_1} once more.

\bigskip\noindent
\emph{Case 3: $p \mid \alpha_{14}$}.\quad
In this case, it is straightforward to compute
\begin{align*}
    \alpha_{14} &= \left( (18 + 36t) - 1 - (7 + 12t) - 1 \right)\left( 1 - (7 + 12t)\right)\\
    &= (9 + 24t)(-6-12t) = -288t^2 - 252t - 54 ,
\end{align*}
hence $\alpha_{41} = 288t^2 + 252t + 54$. Given the fact that
\begin{align*}
    g(\alpha_{41}, n+1-2A_4) &= \frac{1}{2} \left( 36t + 17 + \sqrt{(36t+17)^2 - 4(288t^2 + 252t + 54)}\right)\\
    &= \frac{1}{2} \left( 36t + 17 + \sqrt{144t^2 + 216t + 73}\right)\\
    &= \frac{1}{2} \left( 36t + 17 + \sqrt{(12t+7)^2 + 48t + 24}\right)\\
    &> \frac{1}{2} \left( 36t + 17 + 12t + 7 \right) = 24t + 12
\end{align*}
and
\begin{align*}
    g(\alpha_{41}, n-1) &= \frac{1}{2} \left( 36t + 17 + \sqrt{(36t+17)^2 - 4(288t^2 + 252t + 54)}\right)\\
    &= \frac{1}{2} \left( 36t + 17 + \sqrt{144t^2 + 216t + 73}\right)\\
    &= \frac{1}{2} \left( 36t + 17 + \sqrt{(12t+9)^2 - 8}\right)\\
    &< \frac{1}{2} \left( 36t + 17 + 12t + 9 \right) = 24t + 13 ,
\end{align*}
it quickly follows that if some $p \in \mathbb{N}$ were to satisfy subcondition \textbf{(c)}, we would surely have
\[
    24t + 13 \le p \le 24t + 12 ,
\]
which is obviously impossible.

\bigskip\noindent
\emph{Case 4: $p \mid \alpha_{23}$}.\quad
In this scenario, we get
\begin{align*}
    \alpha_{23} &= \left( (18 + 36t) - 1 - (6 + 12t) - (3 + 12t)\right)\left( (3 + 12t) - (6 + 12t)\right)\\
    &= (8 + 12t)(-3) = -36t - 24,
\end{align*}
together with $\alpha_{32} = 36t + 24$. Suppose that there exists a $p \in \mathbb{N}, \, p \mid \alpha_{23}, \, p > \sqrt{|\alpha_{23}|}$ such that the subconditions \textbf{(a)}, \textbf{(b)} and \textbf{(c)} are all satisfied. Bearing in mind that
\begin{align*}
    g(\alpha_{32}, n + 1 - 2A_3) &= \frac{1}{2} \left( 12t + 13 + \sqrt{(12t+13)^2 - 4(36t+24)}\right)\\
    &= \frac{1}{2} \left( 12t + 13 + \sqrt{144t^2 + 168t + 73}\right)\\
    &= \frac{1}{2} \left( 12t + 13 + \sqrt{(12t + 7)^2 + 24}\right)\\
    &> \frac{1}{2} \left( 12t + 13 + 12t + 7 \right) = 12t + 10
\end{align*}
and
\begin{align*}
    g(\alpha_{32}, n - 1) &= \frac{1}{2} \left( 36t + 17 + \sqrt{(36t + 17)^2 - 4(36t + 24)}\right)\\
    &= \frac{1}{2} \left( 36t + 17 + \sqrt{1296 t^2 + 1080t + 193}\right)\\
    &= \frac{1}{2} \left( 36t + 17 + \sqrt{(36t + 15)^2 - 32}\right)\\
    &< \frac{1}{2} \left( 36t + 17 + 36t + 15 \right) = 36t + 16 ,
\end{align*}
it follows that we certainly have
\begin{equation}\label{xaux_2}
    12t + 11 \le p \le 36t + 15 .
\end{equation}

Now, if $\frac{\alpha_{32}}{p} = 1$, then $p = 36t + 24$, which is impossible due to Eq.~\eqref{xaux_2}. Similarly, if $\frac{\alpha_{32}}{p} = 2$, then $p = 18t + 12$, which leads to a contradiction by virtue of subcondition~\textbf{(a)}. Hence, $\frac{\alpha_{32}}{p} \ge 3$. From here, it follows that $p \le 12t + 8$, thus yielding a contradiction once again due to Eq.~\eqref{xaux_2}.

\bigskip\noindent
\emph{Case 5: $p \mid \alpha_{24}$}.\quad
Here, we have
\begin{align*}
    \alpha_{24} &= \left( (18 + 36t) - 1 - (6 + 12t) - 1 \right)\left( 1 - (6 + 12t)\right)\\
    &= (10 + 24t)(-5-12t) = -288t^2 - 240t - 50 ,
\end{align*}
together with $\alpha_{42} = 288t^2 + 240t + 50$. Since
\begin{align*}
    g(\alpha_{42}, n+1-2A_4) &= \frac{1}{2} \left( 36t + 17 + \sqrt{(36t+17)^2 - 4(288t^2 + 240t + 50)}\right)\\
    &= \frac{1}{2} \left( 36t + 17 + \sqrt{144t^2 + 264t + 89}\right)\\
    &= \frac{1}{2} \left( 36t + 17 + \sqrt{(12t+9)^2 + 48t + 8}\right)\\
    &> \frac{1}{2} \left( 36t + 17 + 12t + 9 \right) = 24t + 13
\end{align*}
and
\begin{align*}
    g(\alpha_{42}, n-1) &= \frac{1}{2} \left( 36t + 17 + \sqrt{(36t+17)^2 - 4(288t^2 + 252t + 54)}\right)\\
    &= \frac{1}{2} \left( 36t + 17 + \sqrt{144t^2 + 264t + 89}\right)\\
    &= \frac{1}{2} \left( 36t + 17 + \sqrt{(12t+11)^2 - 32}\right)\\
    &< \frac{1}{2} \left( 36t + 17 + 12t + 11 \right) = 24t + 14 ,
\end{align*}
we immediately conclude that if some $p \in \mathbb{N}$ were to satisfy subcondition \textbf{(c)}, it would necessarily have to respect the inequality
\[
    24t + 13 \le p \le 24t + 12 ,
\]
which is clearly not possible.

\bigskip\noindent
\emph{Case 6: $p \mid \alpha_{34}$}.\quad
In the final case, we obtain
\begin{align*}
    \alpha_{34} &= \left( (18 + 36t) - 1 - (3 + 12t) - 1 \right)\left( 1 - (3 + 12t)\right)\\
    &= (13 + 24t)(-2-12t) = -288t^2 - 204t - 26,
\end{align*}
as well as $\alpha_{43} = 288t^2 + 204t + 26$. Now, let $p \in \mathbb{N}, \, p > \sqrt{|\alpha_{34}|}$ be a positive divisor of $\alpha_{34}$ which satisfies all the subconditions \textbf{(a)}, \textbf{(b)} and \textbf{(c)}.
Due to
\begin{align*}
    g(\alpha_{43}, n+1-2A_4) &= \frac{1}{2} \left( 36t + 17 + \sqrt{(36t+17)^2 - 4(288t^2 + 204t + 26)}\right)\\
    &= \frac{1}{2} \left( 36t + 17 + \sqrt{144t^2 + 408t + 185}\right)\\
    &= \frac{1}{2} \left( 36t + 17 + \sqrt{(12t+13)^2 + 96t + 16}\right)\\
    &> \frac{1}{2} \left( 36t + 17 + 12t + 13 \right) = 24t + 15
\end{align*}
and
\begin{align*}
    g(\alpha_{43}, n-1) &= \frac{1}{2} \left( 36t + 17 + \sqrt{(36t+17)^2 - 4(288t^2 + 204t + 26)}\right)\\
    &= \frac{1}{2} \left( 36t + 17 + \sqrt{144t^2 + 408t + 185}\right)\\
    &= \frac{1}{2} \left( 36t + 17 + \sqrt{(12t+17)^2 - 104}\right)\\
    &< \frac{1}{2} \left( 36t + 17 + 12t + 17 \right) = 24t + 17 ,
\end{align*}
subcondition \textbf{(c)} dictates that
\[
    24t + 16 \le p \le 24t + 16
\]
must be true. Hence, $p = 24t + 16$. Now, we further have
\[
    -\alpha_{34} = 288t^2 + 204t + 26 = 12t(24t+16) + 12t + 26,
\]
which promptly implies that $p \mid 12t + 26$. If $t = 0$, it is trivial to see that $p \nmid 12t + 26$. Otherwise, if $t \ge 1$, it becomes straightforward to show that $p > 12t + 16 > 0$, which guarantees $p \nmid 12t + 26$ yet again, thus yielding a contradiction. \hfill\qed

\bigskip\noindent
\emph{Proof of Theorem \ref{h_tree_th}, \textbf{(i)}}.\quad
Let $t \in \mathbb{N}_0$ be fixed. By observing that
\begin{alignat*}{5}
    C &= 1 + t, & \qquad A_1 &= 6 + 2t, & \qquad A_2 &= 5 + 2t, & \qquad B_1 &= 2 + t, & \qquad B_2 &= 1,\\
    & & n &= 16 + 6t, & \qquad \frac{n}{2} &= 8 + 3t, & \qquad A_* &= 11 + 4t, & \qquad B_* &= 3 + t ,
\end{alignat*}
it becomes easy to verify that the graph order is even and that conditions  \textbf{(i)} and \textbf{(ii)} from Theorem \ref{main_th_2} certainly hold. Thus, it is enough to demonstrate that conditions \textbf{(iii)}, \textbf{(iv)}, \textbf{(v)} and \textbf{(vi)} are also satisfied by dividing the problem into the eight following cases.

\bigskip\noindent
\emph{Case 1: $p \mid \alpha_{12}$}.\quad
In this case, we have
\begin{align*}
    \alpha_{12} &= \left( (16 + 6t) - 1 - (6 + 2t) - (5 + 2t)\right)\left( (5 + 2t) - (6 + 2t)\right)\\
    &= (4 + 2t)(-1) = -2t - 4,
\end{align*}
as well as $\alpha_{21} = 2t + 4$. Since
\begin{align*}
    g(\alpha_{21}, n + 1 - 2A_2) &= \frac{1}{2} \left( 2t + 7 + \sqrt{(2t+7)^2 - 4(2t+4)}\right)\\
    &= \frac{1}{2} \left( 2t + 7 + \sqrt{4t^2 + 20t + 33}\right)\\
    &= \frac{1}{2} \left( 2t + 7 + \sqrt{(2t+5)^2 + 8}\right)\\
    &> \frac{1}{2} \left( 2t + 7 + 2t + 5 \right) = 2t + 6,
\end{align*}
it follows that if some $p \in \mathbb{N}$ were to satisfy condition \textbf{(iii)}, subcondition \textbf{(c)}, we would then surely have $p > |\alpha_{12} | > 0$. For this reason, such a $p$ could not be a divisor of $\alpha_{12}$.

\bigskip\noindent
\emph{Case 2: $p \mid \beta_{12}$}.\quad
Here, it is easy to see that
\begin{align*}
    \beta_{12} &= \left( (16 + 6t) - 1 - (2 + t) - 1 \right)\left( 1 - (2 + t)\right)\\
    &= (12 + 5t)(-1-t) = -5t^2 - 17t - 12,
\end{align*}
together with $\beta_{21} = 5t^2 + 17t + 12$. Bearing in mind that
\begin{align*}
    g(\beta_{21}, n+1-2B_2) &= \frac{1}{2} \left( 6t + 15 + \sqrt{(6t + 15)^2 - 4(5t^2 + 17t + 12)}\right)\\
    &= \frac{1}{2} \left( 6t + 15 + \sqrt{16t^2 + 112t + 177}\right)\\
    &= \frac{1}{2} \left( 6t + 15 + \sqrt{(4t+13)^2 + 8t + 8}\right)\\
    &> \frac{1}{2} \left( 6t + 15 + 4t + 13 \right) = 5t + 14
\end{align*}
and
\begin{align*}
    g(\beta_{21}, n-1) &= \frac{1}{2} \left( 6t + 15 + \sqrt{(6t + 15)^2 - 4(5t^2 + 17t + 12)}\right)\\
    &= \frac{1}{2} \left( 6t + 15 + \sqrt{16t^2 + 112t + 177}\right)\\
    &= \frac{1}{2} \left( 6t + 15 + \sqrt{(4t + 15)^2 - 8t - 48}\right)\\
    &< \frac{1}{2} \left( 6t + 15 + 4t + 15 \right) = 5t + 15 ,
\end{align*}
we conclude that if some $p \in \mathbb{N}$ were to satisfy condition \textbf{(iv)}, subcondition \textbf{(c)}, we would obtain
\[
    5t + 15 \le p \le 5t + 14,
\]
which is clearly impossible.

\bigskip\noindent
\emph{Case 3: $p \mid \gamma_{1}$}.\quad
In this scenario, we quickly obtain
\begin{align*}
    \gamma_1 &= \left( (16 + 6t) - 1 - (6 + 2t) - (11 + 4t)\right)\left( (11 + 4t) - (6 + 2t)\right)\\
    &= (-2)(5+2t) = -4t - 10 .
\end{align*}
Suppose that there exists a $p \in \mathbb{N}, \, p \mid \gamma_1, \, p > \sqrt{|\gamma_1|}$ such that all three subconditions \textbf{(a)}, \textbf{(b)} and \textbf{(c)} within condition \textbf{(v)} are satisfied. Due to
\begin{align*}
    g(-\gamma_1, 2A_* - n + 3) &= \frac{1}{2} \left( 2t + 9 + \sqrt{(2t + 9)^2 - 4(4t+10)}\right)\\
    &= \frac{1}{2} \left( 2t + 9 + \sqrt{4t^2 + 20t + 41}\right)\\
    &= \frac{1}{2} \left( 2t + 9 + \sqrt{(2t+5)^2 + 16}\right)\\
    &> \frac{1}{2} \left( 2t + 9 + 2t + 5 \right) = 2t + 7
\end{align*}
and
\begin{align*}
    g(-\gamma_1, 2A_* + 2C - n + 1) &= \frac{1}{2} \left( 4t + 9 + \sqrt{(4t + 9)^2 - 4(4t + 10)}\right)\\
    &= \frac{1}{2} \left( 4t + 9 + \sqrt{16t^2 + 56t + 41}\right)\\
    &= \frac{1}{2} \left( 4t + 9 + \sqrt{(4t+7)^2 - 8}\right)\\
    &< \frac{1}{2} \left( 4t + 9 + 4t + 7 \right) = 4t + 8 ,
\end{align*}
subcondition \textbf{(c)} guarantees that
\begin{equation}\label{aaux_1}
    2t + 8 \le p \le 4t + 7 .
\end{equation}
It is clear that $\gamma_1 < 0$. If $\frac{-\gamma_1}{p} = 1$, then $p = 4t + 10$, which is not possible by virtue of Eq.~\eqref{aaux_1}. On the other hand, if $\frac{-\gamma_1}{p} \ge 2$, then $p \le 2t + 5$, which again contradicts Eq.~\eqref{aaux_1}.

\bigskip\noindent
\emph{Case 4: $p \mid \gamma_{2}$}.\quad
Here, it is straightforward to see that
\begin{align*}
    \gamma_2 &= \left( (16 + 6t) - 1 - (5 + 2t) - (11 + 4t)\right)\left( (11 + 4t) - (5 + 2t)\right)\\
    &= (-1)(6+2t) = -2t - 6 .
\end{align*}
Because of
\begin{align*}
    g(-\gamma_2, 2A_* - n + 3) &= \frac{1}{2} \left( 2t + 9 + \sqrt{(2t + 9)^2 - 4(2t+6)}\right)\\
    &= \frac{1}{2} \left( 2t + 9 + \sqrt{4t^2 + 28t + 57}\right)\\
    &= \frac{1}{2} \left( 2t + 9 + \sqrt{(2t+7)^2 + 8}\right)\\
    &> \frac{1}{2} \left( 2t + 9 + 2t + 7 \right) = 2t + 8 ,
\end{align*}
it follows that if some $p \in \mathbb{N}$ were to satisfy condition \textbf{(v)}, subcondition \textbf{(c)}, we would then have $p > |\gamma_2| > 0$. However, such a $p$ could certainly not be a divisor of $\gamma_2$.

\bigskip\noindent
\emph{Case 5: $p \mid \delta_{11}$}.\quad
It is straightforward to compute that
\begin{align*}
    \delta_{11} &= \left( (16 + 6t) - 1 - (6 + 2t) - (2 + t)\right)\left( (2 + t) - (6 + 2t)\right) + (1 + t)(8 + 3t)\\
    &= (7 + 3t)(-4-t) + (3t^2 + 11t + 8) = (-3t^2 - 19t - 28) + (3t^2 + 11t + 8)\\
    &= -8t - 20 .
\end{align*}
Now, suppose that there exists a $p \in \mathbb{N}, \, p \mid \delta_{11}, \, p > \sqrt{|\delta_{11}|}$ such that all three subconditions \textbf{(a)}, \textbf{(b)} and \textbf{(c)} within condition \textbf{(vi)} are satisfied. Bearing in mind that
\begin{align*}
    g(-\delta_{11}, n + 1 - 2B_1) &= \frac{1}{2} \left( 4t + 13 + \sqrt{(4t + 13)^2 - 4(8t + 20)}\right)\\
    &= \frac{1}{2} \left( 4t + 13 + \sqrt{16t^2 + 72t + 89}\right)\\
    &= \frac{1}{2} \left( 4t + 13 + \sqrt{(4t + 9)^2 + 8}\right)\\
    &> \frac{1}{2} \left( 4t + 13 + 4t + 9 \right) = 4t + 11
\end{align*}
and
\begin{align*}
    g(-\delta_{11}, n - 1) &= \frac{1}{2} \left( 6t + 15 + \sqrt{(6t + 15)^2 - 4(8t + 20)}\right)\\
    &= \frac{1}{2} \left( 6t + 15 + \sqrt{36t^2 + 148t + 145}\right)\\
    &= \frac{1}{2} \left( 6t + 15 + \sqrt{(6t+13)^2 - 8t - 24}\right)\\
    &< \frac{1}{2} \left( 6t + 15 + 6t + 13 \right) = 6t + 14 ,
\end{align*}
it becomes evident that subcondition \textbf{(c)} yields
\begin{equation}\label{aaux_2}
    4t + 12 \le p \le 6t + 13 .
\end{equation}
Observe that $\delta_{11} < 0$. If $\frac{-\delta_{11}}{p} = 1$, then $p = 8t + 20$, which is impossible due to Eq.~\eqref{aaux_2}. On the other hand, if $\frac{-\delta_{11}}{p} \ge 2$, this leads us to $p \le 4t + 10$, which is not possible once again by virtue of Eq.~\eqref{aaux_2}.

\bigskip\noindent
\emph{Case 6: $p \mid \delta_{12}$}.\quad
In this case, we have
\begin{align*}
    \delta_{12} &= \left( (16 + 6t) - 1 - (6 + 2t) - 1\right)\left( 1 - (6 + 2t)\right) + (1 + t)(8 + 3t)\\
    &= (8 + 4t)(-5-2t) + (3t^2 + 11t + 8) = (-8t^2 - 36t - 40) + (3t^2 + 11t + 8)\\
    &= -5t^2 - 25t - 32 .
\end{align*}
Since
\begin{align*}
    g(-\delta_{12}, n + 1 - 2B_2) &= \frac{1}{2} \left( 6t + 15 + \sqrt{(6t + 15)^2 - 4(5t^2 + 25t + 32)}\right)\\
    &= \frac{1}{2} \left( 6t + 15 + \sqrt{16t^2 + 80t + 97}\right)\\
    &= \frac{1}{2} \left( 6t + 15 + \sqrt{(4t + 9)^2 + 8t + 16}\right)\\
    &> \frac{1}{2} \left( 6t + 15 + 4t + 9 \right) = 5t + 12
\end{align*}
and
\begin{align*}
    g(-\delta_{12}, n - 1) &= \frac{1}{2} \left( 6t + 15 + \sqrt{(6t + 15)^2 - 4(5t^2 + 25t + 32)}\right)\\
    &= \frac{1}{2} \left( 6t + 15 + \sqrt{16t^2 + 80t + 97}\right)\\
    &= \frac{1}{2} \left( 6t + 15 + \sqrt{(4t+11)^2 - 8t - 24}\right)\\
    &< \frac{1}{2} \left( 6t + 15 + 4t + 11 \right) = 5t + 13 ,
\end{align*}
it follows that if some $p \in \mathbb{N}$ were to satisfy condition \textbf{(vi)}, subcondition \textbf{(c)}, we would have
\[
    5t + 13 \le p \le 5t + 12,
\]
which is obviously impossible.

\bigskip\noindent
\emph{Case 7: $p \mid \delta_{21}$}.\quad
Here, we get
\begin{align*}
    \delta_{21} &= \left( (16 + 6t) - 1 - (5 + 2t) - (2 + t)\right)\left( (2 + t) - (5 + 2t)\right) + (1 + t)(8 + 3t)\\
    &= (8 + 3t)(-3-t) + (3t^2 + 11t + 8) = (-3t^2 - 17t - 24) + (3t^2 + 11t + 8)\\
    &= -6t - 16 .
\end{align*}
Suppose that there exists a $p \in \mathbb{N}, \, p \mid \delta_{21}, \, p > \sqrt{|\delta_{21}|}$ for which all three subconditions \textbf{(a)}, \textbf{(b)} and \textbf{(c)} within condition \textbf{(vi)} are satisfied. Due to
\begin{align*}
    g(-\delta_{21}, n + 1 - 2B_1) &= \frac{1}{2} \left( 4t + 13 + \sqrt{(4t + 13)^2 - 4(6t + 16)}\right)\\
    &= \frac{1}{2} \left( 4t + 13 + \sqrt{16t^2 + 80t + 105}\right)\\
    &= \frac{1}{2} \left( 4t + 13 + \sqrt{(4t + 9)^2 + 8t + 24}\right)\\
    &> \frac{1}{2} \left( 4t + 13 + 4t + 9 \right) = 4t + 11
\end{align*}
and
\begin{align*}
    g(-\delta_{21}, n - 1) &= \frac{1}{2} \left( 6t + 15 + \sqrt{(6t + 15)^2 - 4(6t + 16)}\right)\\
    &= \frac{1}{2} \left( 6t + 15 + \sqrt{36t^2 + 156t + 161}\right)\\
    &= \frac{1}{2} \left( 6t + 15 + \sqrt{(6t+13)^2 - 8}\right)\\
    &< \frac{1}{2} \left( 6t + 15 + 6t + 13 \right) = 6t + 14 ,
\end{align*}
we immediately obtain that subcondition \textbf{(c)} implies
\begin{equation}\label{aaux_3}
    4t + 12 \le p \le 6t + 13 .
\end{equation}
Since $\delta_{21} < 0$, it is clear that if $\frac{-\delta_{21}}{p} = 1$, then $p = 6t + 16$, which is not possible by virtue of Eq.~\eqref{aaux_3}. Similarly, if $\frac{-\delta_{21}}{p} \ge 2$, we directly have $p \le 3t + 8$, which leads us to a contradiction once more due to Eq.~\eqref{aaux_3}.

\bigskip\noindent
\emph{Case 8: $p \mid \delta_{22}$}.\quad
In the final case, it is clear that
\begin{align*}
    \delta_{22} &= \left( (16 + 6t) - 1 - (5 + 2t) - 1\right)\left( 1 - (5 + 2t)\right) + (1 + t)(8 + 3t)\\
    &= (9 + 4t)(-4-2t) + (3t^2 + 11t + 8) = (-8t^2 - 34t - 36) + (3t^2 + 11t + 8)\\
    &= -5t^2 - 23t - 28 .
\end{align*}
Given the fact that
\begin{align*}
    g(-\delta_{22}, n + 1 - 2B_2) &= \frac{1}{2} \left( 6t + 15 + \sqrt{(6t + 15)^2 - 4(5t^2 + 23t + 28)}\right)\\
    &= \frac{1}{2} \left( 6t + 15 + \sqrt{16t^2 + 88t + 113}\right)\\
    &= \frac{1}{2} \left( 6t + 15 + \sqrt{(4t + 9)^2 + 16t + 32}\right)\\
    &> \frac{1}{2} \left( 6t + 15 + 4t + 9 \right) = 5t + 12
\end{align*}
and
\begin{align*}
    g(-\delta_{22}, n - 1) &= \frac{1}{2} \left( 6t + 15 + \sqrt{(6t + 15)^2 - 4(5t^2 + 23t + 28)}\right)\\
    &= \frac{1}{2} \left( 6t + 15 + \sqrt{16t^2 + 88t + 113}\right)\\
    &= \frac{1}{2} \left( 6t + 15 + \sqrt{(4t+11)^2 - 8}\right)\\
    &< \frac{1}{2} \left( 6t + 15 + 4t + 11 \right) = 5t + 13 ,
\end{align*}
we conclude that if some $p \in \mathbb{N}$ were to satisfy condition \textbf{(vi)}, subcondition \textbf{(c)}, we would obtain
\[
    5t + 13 \le p \le 5t + 12,
\]
which is not possible. \hfill\qed

\bigskip
\begin{remark}
    Experimental results suggest that $\mathrm{X}(7 + 33t, 6 + 33t, 3, 1)$ might also be a TI tree for each $t \in \mathbb{N}_0$. Although Theorem \ref{main_th_1} could indeed be applied to potentially verify this statement, the approach disclosed in Section \ref{sc_families} cannot be used in order to carry out the proof.
\end{remark}

\section*{Acknowledgements}
The authors would like to thank Dragan Stevanović, Žana Kovijanić Vukićević and Goran Popivoda for their useful remarks and help in finding all the infinite families of X-trees and H-trees. The authors would also like to express their gratitude to the University of Montenegro for the overall support given throughout the duration of our research.

\section*{Conflict of interest}
The authors declare that they have no conflict of interest.

\appendix

\section{Computed-assisted proof of Theorems \ref{x_tree_th} and \ref{h_tree_th}}\label{python_xh}

\begin{lstlisting}[language = Python, frame = trBL, escapeinside={(*@}{@*)}, aboveskip=10pt, belowskip=10pt, numbers=left, rulecolor=\color{black}]
"""
This Python script provides the computer-assisted proof of Theorems 8 and 9.

Throughout the entire script, Z[t] polynomials will be represented via
one-dimensional integer `np.ndarray` objects. More precisely, a `np.ndarray`
object bearing the form [a_0 a_1 a_2 ... a_{n-1}] shall represent the
polynomial a_0 + a_1 t + a_2 t^2 + ... + a_{n-1} t^{n-1}. The term
coefficients will always appear according to the ascending powers of t.
"""

import numpy as np


def lower_square_root_approximation(input_argument: np.ndarray) -> np.ndarray:
    """
    This function accepts a Z[t] polynomial A(t) = a_0 + a_1 t + a_2 t^2 with
    nonnegative coefficients, and outputs a Z[t] polynomial B(t) = b_0 + b_1 t
    with an even nonnegative b_1 and odd positive b_0, such that the strict
    inequality B(t)^2 < A(t) holds for each nonnegative integer t.

    :arg input_argument: An integer `np.ndarray` object of the shape (3,) that
        represents the input Z[t] polynomial A(t).

    :return: An integer `np.ndarray` object of the shape (2,) that represents
        the output Z[t] polynomial B(t).
    """

    # Check to make sure that the function input is valid.
    assert isinstance(input_argument, np.ndarray)
    assert input_argument.dtype == int and input_argument.shape == (3,)
    assert np.all(input_argument >= 0)

    # Let `leading_coefficient` be the greatest even nonnegative integer whose
    # square does not exceed a_2.
    leading_coefficient = int(np.floor(np.sqrt(input_argument[2])).item())
    if leading_coefficient % 2 == 1:
        leading_coefficient -= 1

    # Let `free_term` be the greatest possible odd nonnegative integer.
    if leading_coefficient == 0:
        # If `leading_coefficient` equals zero, then the square of the free
        # term should not exceed a_0.
        free_term = np.sqrt(input_argument[0])
    else:
        # If `leading_coefficient` is positive, then the square of the free
        # term should not exceed a_0, nor should the expression
        # 2 * `leading_coefficient` * `free_term`
        # exceed a_1, according to the binomial formula.
        free_term = min(
            input_argument[1] / 2.0 / leading_coefficient,
            np.sqrt(input_argument[0]),
        )

    # Make sure that the free term is an odd integer.
    free_term = int(np.floor(free_term).item())
    if free_term % 2 == 0:
        free_term -= 1

    # Raise an error if b_0 and b_1 are not both nonnegative.
    result = np.array([free_term, leading_coefficient], dtype=int)
    assert np.all(result >= 0)

    # Check to make sure that B(t)^2 <= A(t) for any nonnegative integer t.
    result_squared = np.outer(result, result)
    result_squared = np.array(
        [
            result_squared[0, 0],
            result_squared[0, 1] + result_squared[1, 0],
            result_squared[1, 1],
        ],
        dtype=int,
    )
    assert np.all(result_squared <= input_argument)
    # Make sure that the inequality is strict.
    if not result_squared[0] < input_argument[0]:
        result[0] -= 2
        # Remember, b_0 needs to remain nonnegative!
        assert np.all(result >= 0)

    return result


def upper_square_root_approximation(input_argument: np.ndarray) -> np.ndarray:
    """
    This function accepts a Z[t] polynomial A(t) = a_0 + a_1 t + a_2 t^2 with
    nonnegative coefficients, and outputs a Z[t] polynomial B(t) = b_0 + b_1 t
    with an even positive b_1 and odd positive b_0, such that the strict
    inequality B(t)^2 > A(t) holds for each nonnegative integer t.

    :arg input_argument: An integer `np.ndarray` object of the shape (3,) that
        represents the input Z[t] polynomial A(t).

    :return: An integer `np.ndarray` object of the shape (2,) that represents
        the output Z[t] polynomial B(t).
    """

    # Check to make sure that the function input is valid.
    assert isinstance(input_argument, np.ndarray)
    assert input_argument.dtype == int and input_argument.shape == (3,)
    assert np.all(input_argument >= 0)

    # Let `leading_coefficient` be the least even positive integer whose
    # square is at least a_2.
    leading_coefficient = max(
        int(np.ceil(np.sqrt(input_argument[2])).item()), 2
    )
    if leading_coefficient % 2 == 1:
        leading_coefficient += 1

    # Let `free_term` be the lowest possible odd nonnegative integer. The
    # square of the free term should be at least a_0 and the expression
    # 2 * `leading_coefficient` * `free_term`
    # needs to be at least a_1, according to the binomial formula.
    free_term = max(
        input_argument[1] / 2.0 / leading_coefficient,
        np.sqrt(input_argument[0]),
    )
    # Make sure that the free term is an odd positive integer.
    free_term = int(np.ceil(free_term).item())
    if free_term % 2 == 0:
        free_term += 1

    # Raise an error if b_0 and b_1 are not both positive.
    result = np.array([free_term, leading_coefficient], dtype=int)
    assert np.all(result > 0)

    # Check to make sure that B(t)^2 >= A(t) for any nonnegative integer t.
    result_squared = np.outer(result, result)
    result_squared = np.array(
        [
            result_squared[0, 0],
            result_squared[0, 1] + result_squared[1, 0],
            result_squared[1, 1],
        ],
        dtype=int,
    )
    assert np.all(result_squared >= input_argument)
    # Make sure that the inequality is strict.
    if not result_squared[0] > input_argument[0]:
        result[0] += 2

    return result


def lower_g_function_approximation(
    c1: np.ndarray, c2: np.ndarray
) -> np.ndarray:
    """
    This function accepts a Z[t] polynomial C1(t) of degree at most two and a
    Z[t] polynomial C2(t) = c_{20} + c_{21} t, and outputs a Z[t] polynomial
    C3(t) such that C3(t) <= g(C1(t), C2(t)) is true for any nonnegative
    integer t. The coefficient c_{21} is required to be even, while c_{20}
    needs to be odd. Also, the expression C2(t)^2 - 4 * C1(t) is required to
    be a polynomial with nonnegative coefficients.

    :arg c1: An integer `np.ndarray` object of the shape (3,) that represents
        the input Z[t] polynomial C1(t).
    :arg c2: An integer `np.ndarray` object of the shape (2,) that represents
        the input Z[t] polynomial C2(t).

    :return: An integer `np.ndarray` object of the shape (2,) that represents
        the output Z[t] polynomial C3(t).
    """

    # Check to make sure that the function inputs are valid.
    assert isinstance(c1, np.ndarray)
    assert c1.dtype == int and c1.shape == (3,)
    assert isinstance(c2, np.ndarray)
    assert c2.dtype == int and c2.shape == (2,)
    assert c2[0] % 2 == 1 and c2[1] % 2 == 0

    c2_squared = np.outer(c2, c2)
    c2_squared = np.array(
        [
            c2_squared[0, 0],
            c2_squared[0, 1] + c2_squared[1, 0],
            c2_squared[1, 1],
        ],
        dtype=int,
    )
    square_root_argument = c2_squared - 4 * c1

    # Check to make sure that C(t)^2 - 4 * C1(t) has nonnegative coefficients.
    assert np.all(square_root_argument >= 0)
    root_result = lower_square_root_approximation(square_root_argument)

    # Let `result` be the linear polynomial that represents the obtained
    # strict lower bound approximation for g(C1(t), C2(t)).
    result = (c2 + root_result) // 2
    # Since we are dealing exclusively with integers, we may add one to the
    # aforementioned polynomial in order to obtain a (potentially not strict)
    # lower bound approximation for g(C1(t), C2(t)).
    result = result + np.array([1, 0], dtype=int)

    return result


def upper_g_function_approximation(
    c1: np.ndarray, c2: np.ndarray
) -> np.ndarray:
    """
    This function accepts a Z[t] polynomial C1(t) of degree at most two and a
    Z[t] polynomial C2(t) = c_{20} + c_{21} t, and outputs a Z[t] polynomial
    C3(t) such that C3(t) >= g(C1(t), C2(t)) is true for any nonnegative
    integer t. The coefficient c_{21} is required to be even, while c_{20}
    needs to be odd. Also, the expression C2(t)^2 - 4 * C1(t) is required to
    be a polynomial with nonnegative coefficients.

    :arg c1: An integer `np.ndarray` object of the shape (3,) that represents
        the input Z[t] polynomial C1(t).
    :arg c2: An integer `np.ndarray` object of the shape (2,) that represents
        the input Z[t] polynomial C2(t).

    :return: An integer `np.ndarray` object of the shape (2,) that represents
        the output Z[t] polynomial C3(t).
    """

    # Check to make sure that the function inputs are valid.
    assert isinstance(c1, np.ndarray)
    assert c1.dtype == int and c1.shape == (3,)
    assert isinstance(c2, np.ndarray)
    assert c2.dtype == int and c2.shape == (2,)
    assert c2[0] % 2 == 1 and c2[1] % 2 == 0

    c2_squared = np.outer(c2, c2)
    c2_squared = np.array(
        [
            c2_squared[0, 0],
            c2_squared[0, 1] + c2_squared[1, 0],
            c2_squared[1, 1],
        ],
        dtype=int,
    )
    square_root_argument = c2_squared - 4 * c1

    # Check to make sure that C(t)^2 - 4 * C1(t) has nonnegative coefficients.
    assert np.all(square_root_argument >= 0)
    root_result = upper_square_root_approximation(square_root_argument)

    # Let `result` be the linear polynomial that represents the obtained
    # strict upper bound approximation for g(C1(t), C2(t)).
    result = (c2 + root_result) // 2
    # Since we are dealing exclusively with integers, we may subtract one from
    # the aforementioned polynomial in order to obtain a (potentially not
    # strict) upper bound approximation for g(C1(t), C2(t)).
    result = result - np.array([1, 0], dtype=int)

    return result


def confirm_p_value_nonexistence(
    auxiliary_term_1: np.ndarray,
    auxiliary_term_2: np.ndarray,
    auxiliary_term_3: np.ndarray,
    auxiliary_term_4: np.ndarray,
    g_argument_1: np.ndarray,
    g_argument_2: np.ndarray,
    g_argument_3: np.ndarray,
    g_argument_4: np.ndarray,
):
    """
    This function confirms that either condition (iii) from Theorem 6 or one
    of the conditions (iii), (iv), (v) and (vi) from Theorem 7 is satisfied.
    Four input Z[t] polynomials T1(t), T2(t), T3(t) and T4(t) of degree at
    most one are used to construct the fundamental Z[t] polynomial F(t) of
    degree at most two via formula

    F(t) = T1(t) * T2(t) + T3(t) * T4(t).

    The polynomial F(t) signifies one of the values alpha_{ij}, beta_{ij},
    gamma_i and delta_{ij} that appear throughout Theorems 6 and 7. Thus, the
    function confirms that for any nonnegative integer t, there does not exist
    a positive divisor p of F(t) such that p > sqrt(F(t)) and the according
    subconditions (a), (b) and (c) all simultaneously hold. The (b) and (c)
    conditions are represented by the inequalities

    * g(F(t), C1(t)) <= p <= g(F(t), C3(t)),
    * g(-F(t), C2(t)) <= p <= g(-F(t), C4(t)),

    where the according Z[t] polynomials C1(t), C2(t), C3(t) and C4(t) of
    degree at most one are to be prepared as input arguments. If the function
    confirms that all the conditions are satisfied, then no exception is
    thrown. Otherwise, an exception is raised.

    :arg auxiliary_term_1: An integer `np.ndarray` object of the shape (2,)
        that represents the auxiliary input Z[t] polynomial T1(t).
    :arg auxiliary_term_2: An integer `np.ndarray` object of the shape (2,)
        that represents 3he auxiliary input Z[t] polynomial T2(t).
    :arg auxiliary_term_3: An integer `np.ndarray` object of the shape (2,)
        that represents the auxiliary input Z[t] polynomial T3(t).
    :arg auxiliary_term_4: An integer `np.ndarray` object of the shape (2,)
        that represents the auxiliary input Z[t] polynomial T4(t).
    :arg g_argument_1: An integer `np.ndarray` object of the shape (2,) that
        determines the input Z[t] polynomial C1(t).
    :arg g_argument_2: An integer `np.ndarray` object of the shape (2,) that
        determines the input Z[t] polynomial C2(t).
    :arg g_argument_3: An integer `np.ndarray` object of the shape (2,) that
        determines the input Z[t] polynomial C3(t).
    :arg g_argument_4: An integer `np.ndarray` object of the shape (2,) that
        determines the input Z[t] polynomial C4(t).
    """

    # Check to make sure that the function inputs are valid.
    for item in [
        auxiliary_term_1,
        auxiliary_term_2,
        auxiliary_term_3,
        auxiliary_term_4,
        g_argument_1,
        g_argument_2,
        g_argument_3,
        g_argument_4,
    ]:
        assert isinstance(item, np.ndarray)
        assert item.dtype == int and item.shape == (2,)

    # Let `fundament` contain the fundamental polynomial F(t).
    fundament = np.outer(auxiliary_term_1, auxiliary_term_2) + np.outer(
        auxiliary_term_3, auxiliary_term_4
    )
    fundament = np.array(
        [
            fundament[0, 0],
            fundament[0, 1] + fundament[1, 0],
            fundament[1, 1],
        ],
        dtype=int,
    )

    # Print the fundamental polynomial F(t).
    print(f"\tFundamental polynomial: {fundament}.")

    # Use the approximation functions `lower_g_function_approximation` and
    # `upper_g_function_approximation` in order to obtain the lower and upper
    # bounds for the positive divisor p, in accordance with subconditions (b)
    # and (c).
    lower_bound_1 = lower_g_function_approximation(fundament, g_argument_1)
    lower_bound_2 = lower_g_function_approximation(-fundament, g_argument_2)
    upper_bound_1 = upper_g_function_approximation(fundament, g_argument_3)
    upper_bound_2 = upper_g_function_approximation(-fundament, g_argument_4)

    # Let `lower_bound` contain the better of the two lower bounds.
    lower_bound = lower_bound_1
    if lower_bound_2[1] > lower_bound_1[1] or np.all(
        lower_bound_2 >= lower_bound_1
    ):
        lower_bound = lower_bound_2

    # Let `upper_bound` contain the better of the two upper bounds.
    upper_bound = upper_bound_1
    if upper_bound_2[1] < upper_bound_1[1] or np.all(
        upper_bound_2 <= upper_bound_1
    ):
        upper_bound = upper_bound_2

    # Print the two selected bounds.
    print(f"\tLower bound: {lower_bound}.")
    print(f"\tUpper bound: {upper_bound}.")

    # If the lower bound is strictly greater than the upper bound, then it is
    # obvious that no positive divisor p can exist between them.
    if np.all(lower_bound >= upper_bound) and lower_bound[0] > upper_bound[0]:
        return

    # If the fundamental polynomial F(t) is of degree at most one, there
    # exists a convenient way to show that no positive divisor p exists...
    if fundament[2] == 0:
        # Let `absolute_fundament` be the Z[t] polynomial A(t) of degree at
        # most one such that it contains the absolute value of F(t) for each
        # nonnegative integer t. It is checked to make sure that this
        # polynomial has nonnegative coefficients and a positive free term.
        absolute_fundament = fundament[:2]
        if absolute_fundament[0] < 0:
            absolute_fundament = -absolute_fundament
        assert np.all(absolute_fundament >= 0) and absolute_fundament[0] > 0

        # If the lower bound is strictly greater than the value A(t) for each
        # nonnegative integer t, then it is clear that no suitable positive
        # divisor p of A(t) can exist.
        if (
            np.all(lower_bound >= absolute_fundament)
            and lower_bound[0] > absolute_fundament[0]
        ):
            return

        # We verify that the value A(t) is strictly greater than the upper
        # bound. This means that A(t) / p cannot be equal to one.
        assert (
            np.all(absolute_fundament >= upper_bound)
            and absolute_fundament[0] > upper_bound[0]
        )

        # If one half of A(t) is strictly lower than the lower bound, this
        # means that A(t) / p cannot be greater than or equal to two. Thus,
        # there exists no desired positive divisor p.
        if (
            np.all(absolute_fundament <= 2 * lower_bound)
            and absolute_fundament[0] < 2 * lower_bound[0]
        ):
            return

        # Otherwise, we confirm that all the coefficients of A(t) are
        # divisible by four. From here, it follows that if A(t) / p = 2, then
        # p is surely even, yielding an infeasible selection of p due to
        # subcondition (a). For this reason, we can safely assume that
        # A(t) / p is at least three.
        assert np.all(absolute_fundament % 4 == 0)

        # If one third of A(t) is strictly lower than the lower bound, this
        # means that A(t) / p cannot be greater than or equal to three.
        # Therefore, there is no desired positive divisor p.
        if (
            np.all(absolute_fundament <= 3 * lower_bound)
            and absolute_fundament[0] < 3 * lower_bound[0]
        ):
            return

        # Otherwise, we confirm that the free term of A(t) is not divisible by
        # three, while all of its other coefficients are, implying that the
        # value A(t) is surely not divisible by three for any nonnegative
        # integer t. In this case, it is clear that the integer A(t) / p
        # cannot be equal to three. From here, it follows that A(t) / p is
        # necessarily at least four.
        assert (
            absolute_fundament[1] % 3 == 0 and absolute_fundament[0] % 3 != 0
        )

        # Finally, we confirm that one fourth of A(t) is strictly lower than
        # the lower bound, which guarantees that A(t) / p cannot be greater
        # than or equal to four. This implies that there exists no desired
        # positive divisor p of A(t).
        assert (
            np.all(absolute_fundament <= 4 * lower_bound)
            and absolute_fundament[0] < 4 * lower_bound[0]
        )

        return

    # Now, if F(t) is of degree two, we need to have that the lower and upper
    # bound share the same leading coefficient. If they do not, then this
    # approach is simply not capable of confirming the nonexistence of the
    # corresponding positive divisor p.
    assert lower_bound[1] == upper_bound[1]
    leading_coefficient = lower_bound[1]
    # Check to make sure that the leading coefficient is positive.
    assert leading_coefficient > 0
    # Check to make sure that the leading coefficient of the fundamental
    # polynomial F(t) is divisible by the (shared) leading coefficient of the
    # lower and upper bound. If it is not, then once again, we have that the
    # incoming method fails.
    assert fundament[2] % leading_coefficient == 0

    # Since the lower and upper bound share the leading coefficient, this
    # implies that there are only finitely many integer values that a
    # potential feasible positive divisor p could attain. Each of them should
    # be inspected and discarded separately.
    for free_term in range(lower_bound[0], upper_bound[0] + 1):
        # Let `p_expression` contain a linear Z[t] polynomial P(t) which
        # represents a potential expression that the divisor p could equal.
        p_expression = np.array([free_term, leading_coefficient], dtype=int)

        # Let `remainder` be a polynomial yielding a value that is divisible
        # by P(t) if and only if F(t) is, for any nonnegative integer t.
        remainder = fundament.copy()
        # We may assume that the coefficient corresponding to the power two is
        # nonnegative.
        if remainder[2] < 0:
            remainder = -remainder

        # Here, we practically do polynomial division and subtract t * P(t) as
        # many times as needed.
        while remainder[2] >= leading_coefficient:
            remainder[1:] -= p_expression

        # The newly computed coefficient corresponding to the power two must
        # be equal to zero due to the aforementioned assumption made regarding
        # the divisibility of leading coefficients.
        assert remainder[2] == 0
        remainder = remainder[:2]
        # We may now assume that the new coefficient corresponding to the
        # power one is nonnegative.
        if remainder[1] < 0:
            remainder = -remainder

        # We now do another polynomial division step and subtract P(t) as many
        # times as needed.
        while remainder[1] >= leading_coefficient:
            remainder -= p_expression

        # Finally, we obtain a linear remainder polynomial whose coefficient
        # corresponding to the power one must be nonnegative. Let this
        # polynomial be denoted as R(t).
        assert remainder[1] >= 0

        # If the nonnegative integer t equals zero, then the polynomials P(t)
        # and R(t) are simply equal to their free terms, respectively.
        t_value = 0
        p_value = free_term
        r_value = remainder[0]

        # We will finish the confirmation procedure by manually checking
        # finitely many starting values of t up until we have that the
        # condition P(t) > R(t) >= 0 is satisfied. It is obvious that for any
        # nonnegative integer t for which this condition holds, the integer
        # R(t) cannot be divisible by P(t), which implies that neither can
        # F(t). However, bearing in mind that as the t values increase, the
        # P(t) and R(t) values also increase (or stay the same), with P(t)
        # increasing strictly faster than R(t), it becomes sufficient to
        # locate the starting value of t for which the aforementioned
        # condition is indeed satisfied. The remaining values before it have
        # to be inspected manually and for each of them, it should be
        # explicitly demonstrated that no desired positive divisor p exists.
        while not p_value > r_value and r_value >= 0:
            # If the integer R(t) is not divisible by the integer P(t), then
            # there is nothing else to check. But if it is...
            if r_value % p_value == 0:
                # We have to compute the F(t) / P(t) value and compare its
                # parity with that of P(t). If these two integers have the
                # same parity, then the obtained positive divisor p is not
                # feasible, by virtue of subcondition (a). If the integers are
                # of different parities, then our method does not work again.
                fundament_value = np.inner(
                    fundament, [1, t_value, t_value * t_value]
                ).item()
                other_value = fundament_value // p_value
                assert (p_value - other_value) % 2 == 0

            # The t, P(t) and R(t) values should all be updated accordingly.
            t_value += 1
            p_value += leading_coefficient
            r_value += remainder[1]


def confirm_starlike_family(family_specification: np.ndarray):
    """
    This function confirms that all the members of the given family of
    starlike trees are transmission irregular. The function accepts an
    infinite family of starlike trees and implements Theorem 6 in order to
    perform the necessary confirmation. The input family must bear the form

    S(alpha_1 + beta_1 t, alpha_2 + beta_2 t, ..., alpha_k + beta_k t)

    as t ranges over the nonnegative integers, where k >= 3 is a positive
    integer, the alpha integers are all positive and the beta integers are all
    nonnegative. If the function confirms that all the family members are TI,
    then no exception is thrown. An exception is raised if the function is
    unable to make the aforementioned verification.

    :arg family_specification: A two-dimensional integer `np.ndarray` object
        with the shape (2, k) whose elements describe the input alpha and beta
        parameters. The first row contains the alpha_1, alpha_2, ..., alpha_k
        values, respectively, while the second contains the beta_1, beta_2,
        ..., beta_k integers, respectively.
    """

    # Check to make sure that the input family is well specified.
    assert isinstance(family_specification, np.ndarray)
    assert family_specification.dtype == int
    assert (
        len(family_specification.shape) == 2
        and family_specification.shape[0] == 2
        and family_specification.shape[1] >= 3
    )
    assert np.all(family_specification[0]) >= 1
    assert np.all(family_specification[1]) >= 0

    # Print a notification about which starlike tree family is currently being
    # inspected.
    print(f"Inspecting the S-family specified via:\n{family_specification}\n")

    # The order is computed and represented as a linear Z[t] polynomial. It is
    # checked that both of its coefficients are even, thus guaranteeing that
    # the order is even as well.
    order = np.sum(family_specification, axis=1) + np.array([1, 0], dtype=int)
    assert np.all(order % 2 == 0)
    half_order = order // 2

    # The branch lengths are all stored in the `a` variable and represented as
    # linear Z[t] polynomials.
    a = family_specification.T

    # It is checked that the branch lengths are provided in strictly
    # decreasing order. This guarantees that their lengths are all mutually
    # distinct, which implies that condition (i) from Theorem 6 is satisfied.
    for i in range(a.shape[0] - 1):
        assert np.all(a[i] >= a[i + 1]) and a[i, 0] > a[i + 1, 0]

    # It is checked that all the branch lengths are strictly lower than one
    # half of the starlike tree order. This implies that condition (ii) from
    # Theorem 6 holds as well.
    for i in range(a.shape[0]):
        assert np.all(a[i] <= half_order) and a[i, 0] < half_order[0]

    # The auxiliary `confirm_p_value_nonexistence` function is used to verify
    # that condition (iii) from Theorem 6 is satisfied for all the required
    # pairs of integers.
    for i in range(a.shape[0] - 1):
        for j in range(i + 1, a.shape[0]):
            print(f"Resolving case alpha({i + 1}, {j + 1})...")
            confirm_p_value_nonexistence(
                auxiliary_term_1=order
                - np.array([1, 0], dtype=int)
                - a[i]
                - a[j],
                auxiliary_term_2=a[j] - a[i],
                auxiliary_term_3=np.array([0, 0], dtype=int),
                auxiliary_term_4=np.array([0, 0], dtype=int),
                g_argument_1=order + np.array([1, 0], dtype=int) - 2 * a[i],
                g_argument_2=order + np.array([1, 0], dtype=int) - 2 * a[j],
                g_argument_3=order - np.array([1, 0], dtype=int),
                g_argument_4=order - np.array([1, 0], dtype=int),
            )
            print(f"Case alpha({i + 1}, {j + 1}) resolved!\n")

    print("All done!\n")


def confirm_h_family(family_specification: np.ndarray):
    """
    This function confirms that all the members of the given H-family of trees
    are transmission irregular. The function accepts an infinite family of
    H-trees and implements Theorem 7 in order to perform the necessary
    confirmation. The input family must bear the form

    H(alpha_1 + beta_1 t; alpha_2 + beta_2 t, alpha_3 + beta_3 t;
    alpha_4 + beta_4 t, alpha_5 + beta_5 t)

    as t ranges over the nonnegative integers, where the alpha integers are
    all positive and the beta integers are all nonnegative. If the function
    confirms that all the family members are TI, then no exception is thrown.
    An exception is raised if the function is unable to make the
    aforementioned verification.

    :arg family_specification: An integer `np.ndarray` object with the shape
        (2, 5) whose elements describe the input alpha and beta parameters.
        The first row contains the alpha_1, alpha_2, ..., alpha_5 values,
        respectively, while the second contains the beta_1, beta_2, ...,
        beta_5 integers, respectively.
    """

    # Check to make sure that the input family is well specified.
    assert isinstance(family_specification, np.ndarray)
    assert family_specification.dtype == int
    assert family_specification.shape == (2, 5)
    assert np.all(family_specification[0]) >= 1
    assert np.all(family_specification[1]) >= 0

    # Print a notification about which H-tree family is currently being
    # inspected.
    print(f"Inspecting the H-family specified via:\n{family_specification}\n")

    # The order is computed and represented as a linear Z[t] polynomial. It is
    # checked that both of its coefficients are even, thus guaranteeing that
    # the order is even as well.
    order = np.sum(family_specification, axis=1) + np.array([1, 0], dtype=int)
    assert np.all(order % 2 == 0)
    half_order = order // 2

    # The C distance is stored in the `c` variable and represented as a linear
    # Z[t] polynomial. The A_1, A_2 and B_1, B_2 branch lengths are stored in
    # the `a` and `b` variables, respectively, and represented as linear Z[t]
    # polynomials as well.
    c = family_specification[:, 0]
    a = family_specification[:, 1:3].T
    b = family_specification[:, 3:5].T

    # The A_* and B_* values are stored in the `a_star` and `b_star`
    # variables, respectively, as linear Z[t] polynomials.
    a_star = np.sum(a, axis=0)
    b_star = np.sum(b, axis=0)

    # It is checked that A_1 > A_2 and B_1 > B_2, which guarantees that
    # condition (i) from Theorem 7 holds.
    assert np.all(a[0] >= a[1]) and a[0, 0] > a[1, 0]
    assert np.all(b[0] >= b[1]) and b[0, 0] > b[1, 0]

    # We verify that all the A_1, A_2, B_1, B_2 branch lengths are strictly
    # lower than one half of the double starlike tree order.
    for i in range(2):
        assert np.all(a[i] <= half_order) and a[i, 0] < half_order[0]
        assert np.all(b[i] <= half_order) and b[i, 0] < half_order[0]

    # Afterwards, we check that 1 + A_* is strictly greater than one half of
    # the tree order. Thus, condition (ii) from Theorem 7 is necessarily
    # satisfied as well.
    a_star_one = a_star + np.array([1, 0], dtype=int)
    assert np.all(a_star_one >= half_order) and a_star_one[0] > half_order[0]

    # The auxiliary `confirm_p_value_nonexistence` function is used to verify
    # that condition (iii) from Theorem 7 is satisfied.
    print("Resolving case alpha(1, 2)...")
    confirm_p_value_nonexistence(
        auxiliary_term_1=order - np.array([1, 0], dtype=int) - a[0] - a[1],
        auxiliary_term_2=a[1] - a[0],
        auxiliary_term_3=np.array([0, 0], dtype=int),
        auxiliary_term_4=np.array([0, 0], dtype=int),
        g_argument_1=order + np.array([1, 0], dtype=int) - 2 * a[0],
        g_argument_2=order + np.array([1, 0], dtype=int) - 2 * a[1],
        g_argument_3=order - np.array([1, 0], dtype=int),
        g_argument_4=order - np.array([1, 0], dtype=int),
    )
    print("Case alpha(1, 2) resolved!\n")

    # We then apply the same auxiliary function to verify that condition (iv)
    # from Theorem 7 also holds.
    print("Resolving case beta(1, 2)...")
    confirm_p_value_nonexistence(
        auxiliary_term_1=order - np.array([1, 0], dtype=int) - b[0] - b[1],
        auxiliary_term_2=b[1] - b[0],
        auxiliary_term_3=np.array([0, 0], dtype=int),
        auxiliary_term_4=np.array([0, 0], dtype=int),
        g_argument_1=order + np.array([1, 0], dtype=int) - 2 * b[0],
        g_argument_2=order + np.array([1, 0], dtype=int) - 2 * b[1],
        g_argument_3=order - np.array([1, 0], dtype=int),
        g_argument_4=order - np.array([1, 0], dtype=int),
    )
    print("Case beta(1, 2) resolved!\n")

    # Subsequently, we use the same auxiliary function once more in order to
    # check that condition (v) from Theorem 7 is true.
    for i in range(2):
        print(f"Resolving case gamma({i + 1})...")
        confirm_p_value_nonexistence(
            auxiliary_term_1=order
            - np.array([1, 0], dtype=int)
            - a[i]
            - a_star,
            auxiliary_term_2=a_star - a[i],
            auxiliary_term_3=np.array([0, 0], dtype=int),
            auxiliary_term_4=np.array([0, 0], dtype=int),
            g_argument_1=order + np.array([1, 0], dtype=int) - 2 * a[i],
            g_argument_2=2 * a_star - order + np.array([3, 0], dtype=int),
            g_argument_3=order - np.array([1, 0], dtype=int),
            g_argument_4=2 * a_star
            + 2 * c
            - order
            + np.array([1, 0], dtype=int),
        )
        print(f"Case gamma({i + 1}) resolved!\n")

    # Finally, we apply the auxiliary `confirm_p_value_nonexistence` function
    # for the last time in order to verify that condition (vi) from Theorem 7
    # is surely satisfied.
    for i in range(2):
        for j in range(2):
            print(f"Resolving case delta({i + 1}, {j + 1})...")
            confirm_p_value_nonexistence(
                auxiliary_term_1=order
                - np.array([1, 0], dtype=int)
                - a[i]
                - b[j],
                auxiliary_term_2=b[j] - a[i],
                auxiliary_term_3=c,
                auxiliary_term_4=a_star - b_star,
                g_argument_1=order + np.array([1, 0], dtype=int) - 2 * a[i],
                g_argument_2=order + np.array([1, 0], dtype=int) - 2 * b[j],
                g_argument_3=order - np.array([1, 0], dtype=int),
                g_argument_4=order - np.array([1, 0], dtype=int),
            )
            print(f"Case delta({i + 1}, {j + 1}) resolved!\n")

    print("All done!\n")


def confirm_all_x_families():
    """
    This function is used to prove that the four X-tree families

    1. X(7 + 12t, 6 + 12t, 3 + 12t, 1),
    2. X(9 + 12t, 8 + 12t, 5 + 12t, 3),
    3. X(10 + 12t, 9 + 12t, 6 + 12t, 2),
    4. X(11 + 12t, 10 + 12t, 7 + 12t, 1)

    from Theorem 8 all consist of only TI graphs.
    """

    # The families are specified with the help of an integer `np.ndarray`
    # object with the shape (4, 2, 4). Its first dimension corresponds to the
    # four distinct families that are to be inspected, while the other two
    # dimensions align with the family specification mechanism used in the
    # `confirm_starlike_family` function.
    all_family_specifications = np.array(
        [
            [
                [7, 6, 3, 1],
                [12, 12, 12, 0],
            ],
            [
                [9, 8, 5, 3],
                [12, 12, 12, 0],
            ],
            [
                [10, 9, 6, 2],
                [12, 12, 12, 0],
            ],
            [
                [11, 10, 7, 1],
                [12, 12, 12, 0],
            ],
        ],
        dtype=int,
    )

    # The families are inspected one by one with the help of the
    # `confirm_starlike_family` function.
    for family_specification in all_family_specifications:
        confirm_starlike_family(family_specification)


def confirm_all_h_families():
    """
    This function is used to prove that the seven H-tree families

    1. H(1 + t; 6 + 2t, 5 + 2t; 2 + t, 1),
    2. H(2 + t; 7 + 2t, 6 + 2t; 5 + t, 1),
    3. H(1 + 2t; 6 + 4t, 5 + 4t; 1 + 2t, 2),
    4. H(2 + 2t; 8 + 4t, 7 + 4t; 3 + 2t, 1),
    5. H(3 + 2t; 8 + 4t, 7 + 4t; 6 + 2t, 1),
    6. H(1 + t; 7 + 4t, 6 + 4t; 6 + 3t, 1),
    7. H(2; 7 + 4t, 6 + 4t; 5 + 4t, 1)

    from Theorem 9 all consist of only TI graphs. Note that we may safely
    reduce the third family to H(3 + 2t; 10 + 4t, 9 + 4t; 3 + 2t, 2) without
    loss of generality, bearing in mind that its first member coincides with
    the first member of the first family H(1 + t; 6 + 2t, 5 + 2t; 2 + t, 1).
    """

    # The families are specified with the help of an integer `np.ndarray`
    # object with the shape (7, 2, 5). Its first dimension corresponds to the
    # seven distinct families that are to be inspected, while the other two
    # dimensions align with the family specification mechanism used in the
    # `confirm_h_family` function.
    all_family_specifications = np.array(
        [
            [
                [1, 6, 5, 2, 1],
                [1, 2, 2, 1, 0],
            ],
            [
                [2, 7, 6, 5, 1],
                [1, 2, 2, 1, 0],
            ],
            [
                [3, 10, 9, 3, 2],
                [2, 4, 4, 2, 0],
            ],
            [
                [2, 8, 7, 3, 1],
                [2, 4, 4, 2, 0],
            ],
            [
                [3, 8, 7, 6, 1],
                [2, 4, 4, 2, 0],
            ],
            [
                [1, 7, 6, 6, 1],
                [1, 4, 4, 3, 0],
            ],
            [
                [2, 7, 6, 5, 1],
                [0, 4, 4, 4, 0],
            ],
        ],
        dtype=int,
    )

    # The families are inspected one by one with the help of the
    # `confirm_h_family` function.
    for family_specification in all_family_specifications:
        confirm_h_family(family_specification)


if __name__ == "__main__":
    # This is the script entry point.
    confirm_all_x_families()
    confirm_all_h_families()

\end{lstlisting}

\end{document}